\documentclass[11pt]{article}
\usepackage{amsmath,amssymb,latexsym, bbm, times}
\usepackage{amsthm}
\usepackage[T1]{fontenc}
\usepackage[utf8]{inputenc}
\usepackage{rsfso}

\usepackage{amscd}
\usepackage{wrapfig}
\usepackage{tikz}
\usepackage{caption}
\usepackage{subcaption}
\usepackage{pgfplots}
\usepackage{float}
\pgfplotsset{compat=newest}
\usepackage{authblk}
\usepackage{mathrsfs}
\usepackage{bbm}
\usepackage{graphicx}
\usepackage{stmaryrd}
\usepackage{enumitem}
\usepackage[all]{xy}

\usepackage{tocloft}
\usepackage{hyperref}
\hypersetup{
    colorlinks=true,
    linkcolor=blue,
    filecolor=magenta,      
    urlcolor=cyan,
    citecolor=magenta
}
 
\newtheorem{innercustomgeneric}{\customgenericname}
\providecommand{\customgenericname}{}
\newcommand{\newcustomtheorem}[2]{%
  \newenvironment{#1}[1]
  {%
   \renewcommand\customgenericname{#2}%
   \renewcommand\theinnercustomgeneric{##1}%
   \innercustomgeneric
  }
  {\endinnercustomgeneric}
}

\newcustomtheorem{customthm}{Theorem}
\newcustomtheorem{customconj}{Conjecture}  
\newcustomtheorem{customprop}{Proposition}  
\newcustomtheorem{customcor}{Corollary}

\theoremstyle{plain}
\newtheorem{thm}{Theorem}[section]
\newtheorem{lem}[thm]{Lemma}
\newtheorem{prop}[thm]{Proposition}
\newtheorem{cor}[thm]{Corollary}

\newtheorem{thm-def}[thm]{Theorem-Definition}

\theoremstyle{definition}

\theoremstyle{remark}
\newtheorem{rem}{Remark}[section]

\numberwithin{equation}{section}

\newcommand{\gl}{\mathrm{GL}}

\newcommand{\bbg}{\mathbb{G}}

\newcommand{\cf}{\mathcal{F}}
\newcommand{\co}{\mathcal{O}}
\newcommand{\cp}{\mathcal{P}}
\newcommand{\cl}{\mathcal{L}}

\newcommand{\cn}{\mathcal{N}}

\newcommand{\ch}{\mathcal{H}}

\newcommand{\cd}{\mathcal{D}}
\newcommand{\cs}{\mathcal{S}}

\newcommand{\cc}{\mathcal{C}}
\newcommand{\cb}{{\mathcal B}}

\newcommand{\xx}{\mathcal{X}}

\newcommand{\nn}{\mathcal{N}}

\newcommand{\ka}{\mathfrak{a}}
\newcommand{\kg}{\mathfrak{g}}
\newcommand{\kn}{\mathfrak{n}}

\newcommand{\kt}{\mathfrak{t}}

\newcommand{\km}{\mathfrak{m}}

\newcommand{\kl}{\mathfrak{l}}
\newcommand{\kk}{\mathfrak{k}}
\newcommand{\ki}{\mathfrak{i}}
\newcommand{\kr}{\mathfrak{r}}

\newcommand{\bn}{\mathbf{N}}
\newcommand{\bz}{\mathbf{Z}}
\newcommand{\br}{\mathbf{R}}

\newcommand{\bq}{\mathbf{Q}}
\newcommand{\bg}{\mathbf{G}}
\newcommand{\bc}{\mathbf{C}}
\newcommand{\bh}{\mathbf{H}}
\newcommand{\ba}{\mathbf{A}}
\newcommand{\bs}{\mathbf{S}}
\newcommand{\bm}{\mathbf{M}}
\newcommand{\be}{\mathbf{E}}
\newcommand{\bl}{\mathbf{L}}

\newcommand{\bmm}{\mathbf{m}}
\newcommand{\bnn}{\mathbf{n}}

\newcommand{\btt}{\mathbf{t}}

\newcommand{\rv}{\mathrm{v}}
\newcommand{\rw}{\mathrm{w}}
\newcommand{\rc}{\mathrm{c}}
\newcommand{\rd}{\mathrm{d}}
\newcommand{\re}{\mathrm{e}}
\newcommand{\rr}{\mathrm{r}}
\newcommand{\rs}{\mathrm{s}}
\newcommand{\ru}{\mathrm{u}}
\newcommand{\rh}{\mathrm{h}}

\newcommand{\rjj}{\mathrm{J}}
\newcommand{\rdd}{\mathrm{D}}

\newcommand{\val}{\mathrm{val}}

\newcommand{\Ad}{\mathrm{Ad}}
\newcommand{\ad}{\mathrm{ad}}

\newcommand{\diag}{\mathrm{diag}}
\newcommand{\gal}{\mathrm{Gal}}
\newcommand{\Id}{\mathrm{Id}}

\newcommand{\reg}{\mathrm{reg}}

\newcommand{\rss}{\mathrm{rs}}

\newcommand{\Hom}{\mathrm{Hom}}
\newcommand{\ec}{\mathrm{Ec}}

\newcommand{\rk}{\mathrm{rk}}

\newcommand{\ep}{\epsilon}
\newcommand{\vol}{\mathrm{vol}}

\newcommand{\der}{\mathrm{der}}

\newcommand{\fq}{\mathbf{F}_{q}}

\def\ka{{\mathfrak{a}}}

\def\kg{{\mathfrak{g}}}
\def\kh{{\mathfrak{h}}}

\def\kq{{\mathfrak{q}}}
\def\km{{\mathfrak{m}}}
\def\kn{{\mathfrak{n}}}
\def\kz{{\mathfrak{z}}}

\usepackage{xcolor}
\usepackage{version}

\excludeversion{NB}

\newcommand{\supp}{\operatorname{supp}}

\newcommand{\pair}[1]{\langle #1\rangle}

\newcommand{\norm}[1]{|#1|}

\newcommand{\vparen}[1]{\left\{#1\right\}}
\newcommand{\paren}[1]{\left(#1\right)}
\newcommand{\kparen}[1]{\left[#1\right]}

\newcommand{\bsm}{\begin{smallmatrix}}
\newcommand{\esm}{\end{smallmatrix}}

\newcommand{\bfg}{\mathbf{G}}

\def\<{\langle\,}
\def\>{\,\rangle}

\def\kg{\mathfrak g}
\def\kn{\mathfrak n}
\def\<{\langle}
\def\>{\rangle}
\def\={\equiv}
\def\le{\leqslant}
\def\ge{\geqslant}

\def\res{{\rm Res}}

\def\mod{\;{\rm mod}\;}

\def\Hom{{\rm Hom}}

\def\fq{\mathbb{F}_{q}}

\def\rd{{\rm d}}
\def\val{{\rm val}}

\def\gal{{\rm Gal}}

\def\Id{{\rm Id}}

\def\ep{\epsilon}

\author{Zongbin \textsc{Chen}}

\title{\bf On the local behavior of weighted orbital integrals and the affine Springer fibers}

\begin{document}
\maketitle

\begin{abstract}

The paper deals with two problems that seem to be of different nature: the local behavior of the weighted orbital integrals and the dependence of the affine Springer fibers on the root valuation datum of their defining elements. 
The intersection happens when we try to count the number of points on the fundamental domains of the affine Springer fibers: the number of points can be expressed in terms of certain weighted orbital integrals and vice versa. 

For the first problem, we prove the homogeneity of the Arthur-Shalika germ expansion of the weighted orbital integrals,  and establish its parabolic descent. 
This is then applied to our conjecture about the dependence of the Poincar\'e polynomials of the fundamental domains on the root valuation datum. We show that the generating series of the number of points on the fundamental domains with respect to the root valuation datum of the defining elements is a rational fraction, and we bound its denominator with datum associated to the nilpotent orbits.

\end{abstract}

\tableofcontents

\section{Introduction}

Let $k=\fq$ be a finite field, $F=k(\!(\ep)\!)$ the field of Laurent series, $\co=k[\![\ep]\!]$ the ring of integers of $F$.
We fix a separable algebraic closure $\overline{F}$ of $F$.
Let $\val:\overline{F}^{\times}\to \bq$ be the standard valuation on $F$ and let $\norm{{\cdot}}$ be the associated norm. 
Let $\bg$ be a connected reductive algebraic group over $\fq$. Suppose that $\mathrm{char}(k)$ is sufficiently large with respect to the rank of $\bg$.
For a closed algebraic subgroup $\bh$ of $\bg$ over $k$, we denote its Lie algebra by the corresponding gothic letter and we denote $H=\bh(F)$.

Let $\gamma\in \kg$ be a regular semisimple element.  Recall that the \emph{affine Springer fiber} at $\gamma$ is the closed sub ind-$k$-scheme of the \emph{affine grassmannian} $\xx=\bg(\!(\ep)\!)/\bg[\![\ep]\!]$ defined by
$$
\xx_{\gamma}=\big\{g\in \xx\,\mid\,\Ad(g^{-1})\gamma \in \kg[\![\ep]\!]\big\}.
$$ 
The affine Springer fibers are locally finite but not necessarily of finite type over $k$. 
Indeed, let $T$ be the centralizer of $\gamma$ in $G$, then $T$ acts on $\xx_{\gamma}$ by left translation. 
Let $S$ be the maximal $F$-split subtorus of $T$, and let $\Lambda$ be the discrete free Abelian group generated by $\chi_{*}(\ep), \xi\in X_{*}(S)$, then $\Lambda$ acts transitively on the irreducible components of $\xx_{\gamma}$.
In particular, all the irreducible components of $\xx_{\gamma}$ are isomorphic.

Consider the counting points problem on the affine Springer fibers. In our work \cite{chen1}, we have shown that the geometry of $\xx_{\gamma}$ can be reduced to that of its \emph{fundamental domain} $F_{\gamma}$ via the \emph{Arthur-Kottwitz reduction}.
In our work \cite{chen2}, we have explained a transition between counting points on  $F_{\gamma}$ and Arthur's weighted orbital integrals $J_{M_{0}}^{M}(\gamma, \mathbbm{1}_{\km\cap\kk})$. 
The idea is to count the number of points on the \emph{truncated} affine Springer fibers in two ways: by the Arthur-Kottwitz reduction and by the Harder-Narasimhan reduction.  
As a result, we get a recursion involving $|F_{\gamma}^{M}(\fq)|$ and the weighted orbital integral $J_{M_{0}}^{M}(\gamma,\mathbbm{1}_{\km\cap \kk})$, $M\in \cl(M_{0})$. 
The geometry involved is beautiful, but unfortunately the recursion is hard to solve. 
Here we proceed differently, we consider $|F_{\gamma}(\fq)|$ as a weighted count on $\Lambda\backslash\xx_{\gamma}$, the weight factor can be worked out explicitly and related to Arthur's weight factor.

\begin{customthm}{1}[Theorem \ref{f to w}]\label{main f to w}
For any connected reductive algebraic group $\bg$ over $\fq$, we have
$$
|F_{\gamma}(\fq)|=\vol(T^{1})\,|D(\gamma)|^{-\frac{1}{2}}\sum_{\substack{M,\,L\in \cl(M_{0})\\ M\subset L}}
(-1)^{\dim(\ka_{M_{0}}^{M})}J_{M_{0}}^{M}(\gamma, \mathbbm{1}_{\km\cap\kk})\,\rv_{M}^{L}(\ec(x_{0}))\, \re_{L}.
$$
\end{customthm}

The affine Springer fibers behave very nicely with respect to the root valuation datum of their defining elements. For the group $\bg=\gl_{d+1}$, we have proposed a conjecture in this direction: Let $\ba$ be the maximal torus of $\bg$ of the diagonal matrices, let $\mathbf{B}$ be the Borel subgroup of $\bg$ consisting of the upper triangular matrices, let $\{\alpha_{i}\}_{i=1}^{d}$ be the set of simple roots with respect to $\mathbf{B}$. For an integral element $\gamma\in \ka$, it can be shown that up to conjugation by the Weyl group we can put $\gamma$ in the form that for any root $\alpha=\sum_{l=i}^{j}\alpha_{l}$, we have $\val(\alpha(\gamma))=\min\{\alpha_{l}(\gamma)\}_{l=i}^{j}$. The ordered $d$-tuple $(\val(\alpha_{l}(\gamma)))_{l=1}^{d}$ is called the \emph{root valuation datum} of $\gamma$. 

\begin{customconj}{1}[\cite{chen1}]\label{conj rational} 
For $\bfg=\gl_{d+1}$, let $\gamma$ be a diagonal matrix with root valuation datum $\bnn\in \bz_{\ge 0}^{d}$. Then the Poincar\'e polynomial of $F_{\gamma}$ depends only on $\bnn$, we denote it by $P_{\bnn}(t)$. Moreover, the generating series
$$
P(\mathbf{t}):=\sum_{n_{1}=0}^{+\infty}\cdots \sum_{n_{d}=0}^{+\infty}  P_{(n_{1},\cdots,n_{d})}(t) t_{1}^{n_{1}}\cdots t_{d}^{n_{d}}\in \bz [t][\![ t_{1},\cdots,t_{d}]\!]
$$
is a rational fraction, \textup{i.e.} it is an element of $\bz[t](\,t_{1},\cdots,t_{d})$.

\end{customconj}

This is a very strong finiteness statement about the affine Springer fibers. As a consequence, to understand all the affine Springer fibers, it is enough to understand only a finite number of them. The conjecture seems to be out of reach by direct attack, for the lack of effective ways to calculate the Poincar\'e polynomials of $F_{\gamma}$. Nonetheless, for the closely related problem of counting points on $F_{\gamma}$, we are able to establish a result of similar nature:

\begin{customthm}{2}[Theorem \ref{fundamental domain}] \label{main}
For $\bfg=\gl_{d+1}$, let $\gamma$ be a diagonal matrix with root valuation datum $\bnn\in \bz_{\ge 0}^{d}$. 
Then the number of points $|F_{\gamma}(\fq)|$ depends only on $\bnn$, denoted  $F_{\bnn}(q)$. Moreover, the generating series
$$
F(\mathbf{t}):=\sum_{n_{1}=0}^{+\infty}\cdots \sum_{n_{d}=0}^{+\infty}  F_{(n_{1},\cdots,n_{d})}(q) t_{1}^{n_{1}}\cdots t_{d}^{n_{d}}\in \bz [\![ t_{1},\cdots,t_{d}]\!]
$$
is a rational fraction, and we can bound its denominator with datum attached to the nilpotent orbits. 

\end{customthm}

Assuming the purity hypothesis of Goresky, Kottwitz and MacPherson  \cite{gkm1},
which states that $\xx_{\gamma}$ is cohomologically pure in the sense of Grothendieck-Deligne, conjecture \ref{conj rational} will follow from theorem \ref{main}. Indeed, the cohomological purity of $\xx_{\gamma}$ is equivalent to that of $F_{\gamma}$ according to our work \cite{chen1}, and the equivariant cohomology of $F_{\gamma}$ can be calculated as explained in \cite{cl1}, from which we deduce that $P_{\bnn}(t)$ is a polynomial in $t^{2}$ and that $|F_{\gamma}(\fq)|=P_{\bnn}(q^{1/2})$. For more details, see corollary \ref{conj rational reduced}.

The proof of theorem \ref{main} involves the study of the local behavior of weighted orbital integrals. 
By theorem \ref{main f to w}, it is equivalent to a similar rationality assertion for the weighted orbital integrals, which is also  a finiteness statement as we have explained for the generating series $P(q;\mathbf{t})$. Recall that the weighted orbital integrals do satisfy a finiteness property---the \emph{Arthur-Shalika germ expansion}. 
In the group-theoretic setting, Arthur \cite{a1} defined the weighted orbital integrals for the unipotent elements, and showed that the weighted orbital integrals can be expanded in terms of the unipotent ones. In \cite{walds local trace}, Waldspurger adapted Arthur's work to the Lie-algebraic setting. But unfortunately, both authors made no efforts to determine the domain of validity of the expansion. We fill in this gap for the group $\bfg=\gl_{d+1}$. Let $M$ be a Levi subgroup of $G$ containing $\gamma$, let $\cl(M)$ be the set of Levi subgroups of $G$ containing $M$. For $L\in \cl(M)$, let $\nn_{L}$ be the nilpotent cone of $\bl$, and let $[\nn_{L}(F)]$ be a set of representatives for the $L$-conjugacy classes in $\nn_{L}(F)$. 
An element $\gamma\in \kg$ is said to be \emph{integral} if its characteristic polynomial has coefficients in $\co$.
Let $\kg_{0}$ be the set of integral elements in $\kg$.

\begin{customthm}{3}[Theorem \ref{shalika}]

Let $\bfg=\gl_{d+1}$. For each element $u\in \left[\nn_{L}(F)\right], L\in \cl(M)$, there exists unique function $g_{M}^{L}(\,\cdot\,,u)$ on $\km\cap \kg_{0}^{\rm rs}$ such that for any $ f\in \cc_{c}(\kg/\kg(\co))$ and any integral element $\gamma\in \km\cap \kg_{0}^{\rm rs}$, we have
$$
J_{M}(\gamma,f)=\sum_{L\in \cl(M)}\sum_{u\in \left[\nn_{L}(F)\right]}
g_{M}^{L}(\gamma,u)J_{L}(u,f).
$$

\end{customthm}

For the proof, we observe that the weighted orbital integrals $J_{M}(\gamma,f)$ and $J_{L}(u,f)$ have the same behavior under conjugation, to the effect that it can be reduced to the homogeneity of Shalika germ expansion of DeBacker \cite{debacker gl}.  
This expansion will be called \emph{Arthur-Shalika germ expansion}, and the functions $g_{M}^{L}$ will be called \emph{Arthur-Shalika germs}.
Consider the scaling action on $\kg$, which preserves the nilpotent cone and in fact all the nilpotent orbits. By the uniqueness of the Arthur-Shalika germ expansion, we can deduce a homogeneous property of the Arthur-Shalika germs with respect to the scaling action, see theorem \ref{shalika homogeneity}. Moreover, we can consider the germ expansion around an element $s\in \ka_{M}$. This is achieved by a careful analysis of the parabolic descent of the weighted orbital integrals. We refer the reader to theorem \ref{germ descent} for more details.

Coming back to the proof of theorem \ref{main}. As we have explained before, it is equivalent to a rationality statement for a similar generating series of the weighted orbital integrals $J_{A}(\gamma,\mathbbm{1}_{\kg(\co)})$. 
With the Arthur-Shalika germ expansion, this is further reduced to the rationality of a similar generating series for the Arthur-Shalika germs. 
The latter is proved by induction on the root valuation datum of $\gamma$ and on the semisimple rank of $\bg$, using the homogeneity and the parabolic descent of the Arthur-Shalika germs, see theorem \ref{shalika generating series}.

\subsection*{Notations and conventions}

\paragraph*{The structures from the group}

The algebraic subgroups $\bh$ of $\bg$ are always assumed to be defined over $k$, and we use the convention that $H=\bh(F)$. To simplify the notations, we will not refer to the group $\bh$ if the context is clear. This applies particularly to the Levi subgroups and the parabolic subgroups.

Let $M$ be a Levi subgroup of $G$.
We denote by $\cp(M)$ the set of parabolic subgroups of $G$ whose Levi factor is $M$, by $\cl(M)$ the set of Levi subgroups of $G$ containing $M$, and by $\cf(M)$ the set of parabolic subgroups of $G$ containing $M$. For $P\in \cp(M)$, we denote by $P^{-}\in \cp(M)$ the opposite of $P$ with respect to $M$. 
Let 
$$
X^*(M)=\Hom(M, \bbg_m), \quad X_{*}(M)=\Hom(X^{*}(M), \bz).
$$
Let $\ka_M^{*}=X^*(M)\otimes\br$ and $\ka_M=X_{*}(M)\otimes\br$. Let $A_{M}$ be the central  connected component of the center of $M$, then the restriction $X^{*}(M)\to X^{*}(A_{M})$ induces an isomorphism $\ka_{M}^{*}\cong X^{*}(A_{M})\otimes \br$. 
For a Levi subgroup $L\in \cl(M)$, the restriction $X^{*}(L)\to X^{*}(M)$ induces an injection $\ka_{L}^{*}\hookrightarrow \ka_{M}^{*}$, hence a dual surjection $\ka_{M}\twoheadrightarrow \ka_{L}$. On the other hand, the inclusion $A_{L}\subset A_{M}$ induces an inclusion $X_{*}(A_{L})\otimes \br \hookrightarrow X_{*}(A_{M})\otimes \br$ and a dual surjection $X^{*}(A_{M})\otimes \br\twoheadrightarrow X^{*}(A_{L})\otimes \br$. Taken together, we get splittings
$$
\ka_{M}=\ka_{M}^{L} \oplus  \ka_{L},\quad (\ka_{M})^{*}=(\ka_{M}^{L})^{*} \oplus  (\ka_{L})^{*}.
$$
Let $\pi^{L},\,\pi_{L}$ be the projections to the two factors for these two splittings.

To define Arthur's weight factor $\rv_{M}(g)$, we need to choose a Lebesgue measure on $\ka_{M}^{G}$. We fix a $W$-invariant positive definite inner product $\langle\, \cdot\,,\cdot\,\rangle$ on the vector space $\ka_{A}^{G}$, where $W$ is the Weyl group of $G$ with respect to $A$. Notice that $\ka_{A}^{M}$ and $\ka_{M}$ are orthogonal to each other with respect to the inner product for any $M\in \cl(A)$. We fix a Lebesgue measure on $\ka_{M}^{G}$ normalised by the condition that the lattice generated by the orthonormal bases in $\ka_{M}^{G}$ has covolume $1$.

For $P\in \cf(A)$, let $M_{P}$ be its Levi factor containing $A$ and let $N_{P}$ be its unipotent radical. Let $\Phi_{P}$ be the set of roots of $(P,A_{P})$ for  $A_{P}:=A_{M_{P}}$, let $\overline{\Phi}_{P}$ be the reduced root system and let $\Delta_{P}$ be the set of simple roots in $\overline{\Phi}_{P}$, they can be seen as elements in $(\ka_{M}^{G})^{*}$. Let $\Delta_{P}^{\vee}$ be the associated simple coroots. 


\paragraph*{$\mathbf{(G,M)}$-family}

We need Arthur's notion of \emph{$(G,M)$-family}. The reader can find a beautiful exposition in \cite{a2}, \S17. It is a family of smooth functions $(\rr_{P}(\lambda))_{P\in \cp(M)}$ on $i\ka_{M}^{*}\subset \ka_{M}^{*}\otimes \bc$ which satisfy for any adjacent parabolic subgroups $(P,P')$ the property that $\rr_{P}(\lambda)=\rr_{P'}(\lambda)$ for any $\lambda$ on the hyperplane defined by the unique coroot $\beta_{P,P'}^{\vee}$ in $\Delta_{P}^{\vee}\cap (-\Delta_{P'}^{\vee})$. An important source of example comes from the $(G,M)$-orthogonal sets. 
A collection $r=(r_{P})_{P\in \cp(M)}$ with $r_{P}\in \ka_{M}$ is said to be a \emph{$(G,M)$-orthogonal set} if for any adjacent parabolic subgroups $P,P'\in \cp(M)$,
$$
r_{P}-r_{P'}= c_{P,P'}\cdot \beta_{P,P'}^{\vee},
$$
for some $c_{P,P'}\in \br$. It is said to be \emph{positive} if all the $c_{P,P'}$ are positive. Notice that the family of function $\rr=(\rr_{P}(\lambda))_{P\in \cp(M)}$ with 
$$
\rr_{P}(\lambda):=e^{\pair{\lambda,\,r_{P}}},\quad \lambda\in i\ka_{M}^{*},
$$
forms a $(G,M)$-family, called the $(G,M)$-family associated to the $(G,M)$-orthogonal set.

For $Q=LN_{Q}\in \cf(M)$, let
$$
\rr_{R}^{Q}(\lambda)=\rr_{RN_{Q}}(\lambda), \quad \forall \, R\in \cp^{L}(M).
$$ 
It is easy to see that $\big(\rr_{R}^{Q}(\lambda)\big)_{R\in \cp^{L}(M)}$ form a $(L,M)$-family. We call it the \emph{$Q$-facet} of $\rr$, denoted by $\rr^{Q}$. Notice that it can be viewed as a $(G,M)$-family, denoted $\iota(\rr^{Q})$, by $\iota(\rr^{Q})_{P}(\lambda)=\rr_{P\cap L}^{Q}(\lambda)$ for $P\in \cp(M)$.  
On the other hand, let $P\in \cp(M)$ be a parabolic subgroup contained in $Q$, it is easy to see that the restriction of the function $\rr_{P}(\lambda)$ to $\ka_{L}^{*}$ doesn't depend on $P$, we denote it by $\rr_{Q}(\lambda)$. It is easy to see that the family $(\rr_{Q}(\lambda))_{Q\in \cp(L)}$ forms a $(G,L)$-family, we denote it by $\pi_{L}(\rr)$. Similar construction works for the $(G,M)$-orthogonal sets, we preserve the same notations for them. Notice that there is no ambiguity about the notation $\pi_{L}(r)$ because it is also the image of $r$ under the projection $\pi_{L}:\ka_{M}\to \ka_{L}$.

In case that the $(G,M)$-orthogonal set $r$ is positive, we denote by $\be(r)$ the convex hull of the points $r_{P}$, it is a convex polytope in $\ka_{M}$. Its volume with respect to the Lebesgue measure on $\ka_{M}$ can be calculated to be 
\begin{equation*}
\rr_{M}=\lim_{\lambda\to 0}\sum_{P\in \cp(M)} \rd_{P}(\lambda)^{-1} \, \rr_{P}(\lambda),
\end{equation*}
with $\rd_{P}(\lambda)=\vol(\ka_{M}^{G}/\bz(\Delta_{P}^{\vee}))^{-1}\prod_{\alpha\in \Delta_{P}}\lambda(\alpha^{\vee})$. For a general $(G,M)$-family $(\rr_{P}(\lambda))_{P\in \cp(M)}$, we define
$$
\rr_{M}(\lambda)=\sum_{P\in \cp(M)} \rd_{P}(\lambda)^{-1} \, \rr_{P}(\lambda),
$$
for generic $\lambda\in i\ka_{M}^{*}$. It can be shown that the function extends smoothly over all $i\ka_{M}^{*}$. Let
\begin{equation}\label{def gm volume}
\rr_{M}=\lim_{\lambda\to 0} \rr_{M}(\lambda)=\lim_{\lambda\to 0}
\sum_{P\in \cp(M)} \rd_{P}(\lambda)^{-1} \, \rr_{P}(\lambda),
\end{equation}
and we call it the \emph{volume} of the $(G,M)$-family $(\rr_{P}(\lambda))_{P\in \cp(M)}$.

Let $c$ and $d$ be two $(G,M)$-orthogonal sets, we define their sum $c+d$ to be 
$$
(c+d)_{P}=c_{P}+d_{P},\quad \forall P\in \cp(M).
$$
Obviously the sum remains a $(G,M)$-orthogonal set. Similarly, we can define the product of two $(G,M)$-families. If both $c$ and $d$ are positive, and that $0\in \ec(d)$, we have inclusion of polytopes $\be(c)\subset \be(c+d)$, and the complement $\be(c+d)\backslash \be(c)$ can be partitioned as: For $Q=LN_{Q}\in \cf(M)$, let 
$$
\ka^{+}_{Q}=\vparen{a\in \ka_{M}\mid \alpha(a)> 0, \forall \alpha\in \Delta_{Q};\, \beta(a)=0,\forall \beta\in \Phi(L,A_{M})}. 
$$ 
Then we have disjoint partition 
$$
\ka_{M}=\bigsqcup_{Q\in \cf(M)}\ka^{+}_{Q}.
$$ 
Let $\be^{Q}(c)$ be the $Q$-facet of the polytope $\be(c)$, let $R_{Q}(c)=\be^{Q}(c)+\ka^{+}_{Q}$, then 
$$
\ka_{M}=\bigsqcup_{Q\in \cf(M)}R_{Q}(c).
$$
Let $R_{Q}(c,d)=R_{Q}(c)\cap \be(c+d)$, we get a decomposition
$$
\be(c+d)=\bigsqcup_{Q\in \cf(M)}R_{Q}(c,d),
$$
with $R_{G}(c,d)=\be(c)$. The decomposition enables us to express the volume of $\be(c+d)$ in terms of that of $c$ and $d$.

The construction generalizes to the $(G,M)$-families.
Let $\rr,\rs$ be two $(G,M)$-families. Then the collection $\vparen{\rr_{P}(\lambda)\rs_{P}(\lambda)}_{P\in \cp(M)}$ forms a $(G,M)$-family. We denote it by $\rr\cdot \rs$ and call it the product of $\rr$ and $\rs$. From the $(G,M)$-family $\{\rs_{P}(\lambda)\}_{P\in \cp(M)}$, Arthur has defined a smooth function $\rs_{Q}'(\lambda)$ on $\ka_{Q}^{*}$ for $Q\in \cf(M)$, which is too complicated to recall here and we simply refer the reader to \cite{a2}, equation (17.7). 
Let $\rs_{Q}'=\rs_{Q}'(0)$.
Then Arthur's formula for the volume of the product of two $(G,M)$-families states 
\begin{equation}\label{arthur prod}
(\rr\cdot \rs)_{M}=\sum_{Q\in \cf(M)}\rr_{M}^{Q}\cdot\rs_{Q}'.
\end{equation}
Setting $\rr=1$, we get 
\begin{equation}\label{vol der}
\rs_{M}=\sum_{P\in \cp(M)}\rs_{P}'.
\end{equation}
In case that the $(G,M)$-family $\rr$ satisfies the extra conditon
$$
\rr_{M}^{L}=\rr_{M}^{Q},\quad \forall\, L\in \cl(M),Q\in \cp(L), 
$$
the above product formula further simplifies to
\begin{equation}\label{arthur prod simp}
(\rr\cdot \rs)_{M}=\sum_{L\in \cl(M)} \rr_{M}^{L}\cdot\rs_{L}.
\end{equation}

For $L\in \cl(M)$, we will need a descent formula for the volume $\rr_{L}$ of the $(G,L)$-family $(\rr_{R}(\lambda))_{R\in \cp(L)}$. For $L'\in \cl(M)$, we define $\theta_{M}^{G}(L,L')$ to be 0 unless the natural map
$$
\ka_{M}^{L}\oplus \ka_{M}^{L'}\to \ka_{M}^{G}
$$
is an isomorphism, in which case $\theta_{M}^{G}(L,L')$ is the factor by which the product  Haar measure on $\ka_{M}^{L}\oplus \ka_{M}^{L'}$ must be multiplied in order to be equal to the Haar measure on $\ka_{M}^{G}$. Let $\xi$ be a fixed generic vector in $\ka_{M}^{L}$, then for any $L'\in \cl(M)$ satisfying $\theta_{M}^{G}(L,L')\neq 0$, the affine space $\xi+\ka_{M}^{G}$ intersects non-trivially with $\ka_{L'}^{G}\subset \ka_{M}^{G}$ at one point, which belongs to a chamber $\ka_{Q}^{+}$ in $\ka_{L'}^{G}$, for a unique $Q\in \cp(L')$. We will denote it by $Q_{L'}^{\xi}$. Then
\begin{equation}\label{proj vol}
\rr_{L}=\sum_{L'\in \cl(M)} \theta_{M}^{G}(L,L')\cdot \rr_{M}^{Q_{L'}^{\xi}}.
\end{equation}

\paragraph*{Measures on the groups}

To define the weighted orbital integrals in a uniform way, we need to make a uniform choice of Haar measure on $G$ and its closed subgroups. Let $dx$ be the Haar measure on $F$ with $dx(\co)=1$. 
Let $\bh$ be a connected closed algebraic subgroup of $\bg$ over $F$. Following Weil \cite{weil integration}, let $\omega_{\kh}\in \bigwedge^{\rm top}\kh$ be a volume form on $\kh$, then it defines a Haar measure $|\omega_{\kh}|$ on $\kh$.
By left translation, $\omega_{\kh}$ defines a left-invariant volume form $\omega_{\bh}$ on $\bh$. This determines a left-invariant measure $|\omega_{\bh}|$ on $H$. 
We normalize $|\omega_{\kh}|$ and $|\omega_{\bh}|$ by taking a volume form $\omega_{\kh}$ such that $\vol_{|\omega_{\bh}|}(H\cap \bg(\co))=1$.
We denote the resulting measure by $dh$ if the context is clear. 
Given an element $\gamma\in \kg$, it is known that its stabilizer $G_{\gamma}$ for the adjoint action of $G$ is unimodular, hence $G_{\gamma}$ carries a normalized Haar measure $dh$ and we can define an invariant measure $\frac{dg}{dh}$ on the $G$-orbit of $\gamma$.

\paragraph*{Moy-Prasad subgroups}

Let $\cb(G)$ be the Bruhat-Tits building of $G$. For any point $x\in \cb(G)$, Moy and Prasad \cite{mp} defines a decreasing exhaustive filtration $\{\kg_{x,r}\}_{r\in \br}$ of $\kg$. For $r\ge 0$, $\kg_{x,r}$ is an ideal of $\kg$, let $G_{x,r}$ be the associated Lie group, it is a normal subgroup of $G$. 
Let $\kg_{x,r+}=\bigcup_{t>r}\kg_{x,t}$ and let $\kg_{x}(r)=\kg_{x,r}/\kg_{x,r+}$.
For $x=0$, we simplify $\kk=\kg_{0,0}$ and $\kk_{r}=\kg_{0,r}$, and accordingly $K=G_{0,0}, K_{r}=G_{0,r}$ for $r\ge 0$. 
For $r=0$, we abbreviate $\kg_{x}=\kg_{x,0}$ and accordingly $G_{x}=G_{x,0}$.
Let
$$
\kg_{r}=\bigcup_{x\in \cb(G)}\kg_{x,r},\quad r\in \br.
$$
Let $\cc_{c}^{\infty}(\kg)$ be the space of compactly supported locally constant functions on $\kg$. We define its subspaces
$$
\ch_{r}=\sum_{x\in \cb(G)} \cc_{c}(\kg/\kg_{x,r}).
$$
Similar definitions for $\kg_{r+}$ and $\ch_{r+}$.

\paragraph*{Connected components of the affine grassmannian}

We identify $X_{*}(A)$ with $A(F)/A(\co)$ by sending $\chi$ to $\chi(\ep)$. With this identification, the canonical surjection $A(F)\to A(F)/A(\co)$ can be viewed as
\begin{equation}\label{indexT}
A(F)\to X_{*}(A).
\end{equation}

We use $\Lambda_{G}$ to denote the quotient of $X_{*}(A)$ by the coroot lattice of $G$ (the subgroup of $X_{*}(A)$ generated by the coroots of $A$ in $G$). According to Kottwitz \cite{k1}, we have a canonical homomorphism
\begin{equation}\label{indexM}
\nu_{G}: G\to \Lambda_{G},
\end{equation}
which is characterized by the following properties: it is trivial on the image of $G_{\mathrm{sc}}$ in $G$ ($\bg_{\mathrm{sc}}$ is the simply connected cover of $\bg^{\rm der}$), and its restriction to $A$ coincides with the composition of (\ref{indexT}) with the projection of $X_{*}(A)$ to $\Lambda_{G}$. Since the morphism (\ref{indexM}) is trivial on $K$, it descends to a map
$$
\nu_{G}:\xx\to \Lambda_{G},
$$
whose fibers are the connected components of $\xx$. For $\mu\in \Lambda_{G}$, we denote the connected component $\nu_{G}^{-1}(\mu)$ by $\xx^{\mu}$.

\section{The weighted orbital integrals}

We make a brief recall on the weighted orbital integrals for both the group $G$ and its Lie algebra. 
We can work in a broader context: Let $F$ be a complete discrete valued field, let $\bg$ be a connected reductive algebraic group over $F$. Assume that ${\rm char}(F)$ is either $0$ or large enough with respect to $\rk(\bg)$. We fix a minimal Levi subgroup $\bm_{\rm min}\subset \bg$, and a special vertex in the apartment associated to $A_{M_{\rm min}}$ in the Bruhat-Tits building $\cb(G)$ of $G$, let $K$ be the stabilizer in $G$ of this special vertex. We work exclusively with the Levi subgroups $M\in \cl(M_{\rm min})$.

For $P=MN\in \cp(M)$, the Iwasawa decomposition $G=NMK$ defines a parabolic reduction  on the affine grassmannian: 
$$
f_{P}:\xx\to \xx^{M},\quad [g]=[nmk]\mapsto [m],
$$
Let $H_{M}: M\to \ka_{M}$ be the unique group homomorphism\footnote{Our definition differs from the conventional one by a minus sign.} satisfying
$$
\chi(H_{M}(m))=\val(\chi(m)),\quad \forall\, \chi\in X^{*}(M),\,m\in M.
$$
Notice that it is invariant under the right $K$-action, so it induces a map from $\xx^{M}$ to $\ka_{M}$, still denoted by $H_{M}$. Let $H_{P}:\xx\to \ka_{M}$ be the composition 
$$
H_{P}:\xx\xrightarrow{f_{P}}\xx^{M}\xrightarrow{H_{M}}\ka_{M}.
$$
It has been shown by Arthur \cite{a} that for any point $x\in \xx$, the reductions 
$$
\ec_{M}(x):=\big(H_{P}(x)\big)_{P\in \cp(M)}
$$ 
forms a positive $(G,M)$-orthogonal set. Let 
$$
\rv_{P}(\lambda, x)=e^{\pair{\lambda,\, H_{P}(x)}}, \quad \lambda\in i\ka_{M}^{*},
$$
then $\big(\rv_{P}(\lambda, x)\big)_{P\in \cp(M)}$ forms a $(G,M)$-family. Let 
$\rv_{M}(x)$ be its volume as defined by the equation (\ref{def gm volume}), we call it \emph{Arthur's weight factor}. 

Let $x\in M$ be an element satisfying $G_{x}=M_{x}$. For $\phi\in \cc_{c}^{\infty}(G)$, the \emph{weighted orbital integral} $\rjj_{M}(x,\phi)$ is defined as 
\begin{equation}\label{def w int group}
\rjj_{M}(x,\,\phi)=
|\rdd(x)|^{1/2}\int_{G_{x}\backslash G} \phi\big(\Ad(g^{-1})x \big)\rv_{M}(g)\frac{dg}{dt},
\end{equation}
which is well-defined because $\rv_{M}$ is left $M$-invariant and $G_{x}=M_{x}$ is contained in $M$. Here ${\rdd(x)}$ is defined by
$$
\rdd(x):=\det\big(\Ad(x_{s})-\Id\mid \kg/\kg_{x_{s}}\big),
$$
where $x_{s}$ is the semisimple part of $x$.

Correspondingly, for an element $\gamma\in \km$ satisfying $G_{\gamma}=M_{\gamma}$ and a test function $\varphi\in \cc_{c}^{\infty}(\kg)$, the \emph{weighted orbital integral} $J_{M}(\gamma,\varphi)$ is defined as 
\begin{equation}\label{def w int 1}
J_{M}(\gamma,\,\varphi)=
|D(\gamma)|^{1/2}\int_{G_{\gamma}\backslash G} \varphi\big(\Ad(g^{-1})\gamma \big)\rv_{M}(g)\frac{dg}{dt}.
\end{equation}
Again it is well-defined because $\rv_{M}$ is left $M$-invariant and $G_{\gamma}=M_{\gamma}$ is contained in $M$. The function ${D(\gamma)}$ is defined as 
$$
D(\gamma):=\det\big(\ad(\gamma_{s})\mid \kg/\kg_{\gamma_{s}}\big),
$$
where $\gamma_{s}$ is the semisimple part of $\gamma$.

\begin{rem}\label{woi group to lie}

Following Waldspurger \cite{walds local trace}, we can reduce questions about weighted orbital integrals on the Lie algebra to weighted orbital integrals on the group with the following trick: Let $V_{\kg}$ (resp. $V_{G}$) be sufficiently small $G$-invariant neighbourhood of $0\in \kg$ (resp. $e\in G$) such that the exponential map defines a homeomorphism $\exp:V_{\kg}\to V_{G}$. Take $n\in \bn$ sufficiently large such that $\ep^{n}\gamma\in V_{\kg}$ and the function $\varphi_{n}$ on $\kg$ defined by $\varphi_{n}(X)=\varphi(\ep^{-n}X)$ is supported on $V_{\kg}$, let $\phi_{n}$ be the push-forward of $\varphi_{n}$ via the exponential map. Let $x_{n}=\exp(\ep^{n}\gamma)\in M$, comparing (\ref{def w int group}) and (\ref{def w int 1}), we get
$$
J_{M}(\gamma,\,\varphi)=
|\ep^{n\rk(\kg/\kg_{\gamma_{s}})}|^{-1/2}\rjj_{M}(x_{n}, \phi_{n}).
$$
It is often more convenient to make change of variable on the groups. 

\end{rem}

For general element $x\in M$, without the assumption that $G_{x}=M_{x}$, following Arthur \cite{a1}, we perturb $x$ by an element $a\in A_{M}^{\reg}$ sufficiently close to $1$ such that $G_{ax}=M_{ax}$, and define $\rjj_{M}(x,\phi)$ by a limit process.

\begin{thm}[Arthur \cite{a1}]\label{woi unipotent}

For $x\in M$, let $x=x_{s}x_{u}$ be the Jordan decomposition of $x$. For any test function $\phi\in \cc_{c}^{\infty}(G)$, the limit
\begin{equation*}
\lim_{\substack{a\in A_{M}^{\reg}\\ a\to 1}}
\sum_{L\in \cl(M)} \rr_{M}^{L}(x_{u},a)\cdot \rjj_{L}(ax,\phi)
\end{equation*}
exits, where $\rr_{M}^{L}(x_{u},a)$ is the volume of an $(L,M)$-family.
The limit is defined as the weighted orbital integral associated to $x$ and $\phi$, denoted $\rjj_{M}(x,\phi)$. 

\end{thm}

With remark \ref{woi group to lie}, the theorem can be translated to the Lie algebraic setting. For any element $\gamma\in \km$ and test function $\varphi\in \cc_{c}^{\infty}(\kg)$, we define $J_{M}(\gamma,\varphi)$ as the limit
\begin{equation}\label{def w int 2}
J_{M}(\gamma,\varphi):=\lim_{\substack{a\in \ka_{M}^{\reg}\\ a\to 0}}
\sum_{L\in \cl(M)} \rr_{M}^{L}(\exp(\gamma_{u}),\exp(a))\cdot J_{L}(\gamma+a,\phi).
\end{equation}
We simplify the notation $\rr_{M}^{L}(\exp(\gamma_{u}),\exp(a))$ to $\rr_{M}^{L}(\gamma_{u},a)$ or even $\rr_{M}^{L}(\gamma,a)$ if there is no risk of confusion.

We want to emphasize the approximation of $\gamma$ by $\gamma+a$ with $a\in \ka_{M}^{\reg}$  in the definition (\ref{def w int 2}). This is to ensure that the $M$-orbit of $\gamma+a$ approaches the $M$-orbit of $\gamma$. If we used the same approximation to calculate $J_{L}(\gamma,\varphi)$ for $L\in \cl(M)$, we will arrive at $J_{L}(\gamma_{M}^{L},\varphi)$ with $\gamma_{M}^{L}$ the Lusztig-Spaltenstein induction of $\gamma$, because the $L$-orbit of $\gamma+a$  will approach the $L$-orbit of $\gamma_{M}^{L}$ rather than that of $\gamma$. Recall that the \emph{Lusztig-Spaltenstein induction} $\gamma_{M}^{G}$ of an element $\gamma\in \km$ is the union of the $G$-orbits $\{\gamma_{i}\}_{G}$ which intersects the coset $\gamma+\kn$ in a dense open subset, where $\kn$ is the Lie algebra of the unipotent radical of a parabolic subgroup $P\in \cp(M)$. There are only finitely many such orbits and they belong to the same geometric orbit. They can be characterized analytically by the formula
$$
J_{G}(\gamma_{M}^{G},\varphi):=\sum_{i}J_{G}(\gamma_{i},\varphi)=\lim_{\substack{a\in \ka_{M}^{\reg}\\ a\to 0}} J_{G}(\gamma+a,\varphi), \quad \forall\,\varphi\in \cc_{c}^{\infty}(\kg),
$$
which is a reflection of the geometric fact that the $G$-orbit of $\gamma+a$ approaches the $G$-orbits in $\gamma_{M}^{G}$ when $a\in \ka_{M}^{\reg}\mapsto 0$. In general, let $\gamma_{M}^{L}=\bigsqcup_{i}\vparen{\gamma'_{i}}_{L}$ be the decomposition into $L$-orbits, we define
$$
J_{L}(\gamma_{M}^{L},\varphi)=\sum_{i}J_{L}(\gamma'_{i},\varphi).
$$
This is well-defined because $J_{L}(\gamma'_{i},\varphi)$ depends only on the $L$-orbit of $\gamma_{i}'$. Then following the same argument as in proof of \cite{a1}, cor. 6.3,  we can show
\begin{equation}\label{refined appr}
\lim_{\substack{a\in \ka_{M}^{\reg}\\ a\to 0}}
\sum_{L'\in \cl(L)} \rr_{L}^{L'}(\gamma_{u},a)\cdot J_{L'}(\gamma+a,\varphi)=J_{L}(\gamma_{M}^{L},\varphi).
\end{equation}

\section{Homogeneity of the Arthur-Shalika germ expansion}

The great advantage of definition (\ref{def w int 2}) is that $J_{M}(\gamma,\varphi)$ has the same behavior under conjugation, whether $\gamma \in \km$ satisfies the condition $M_{\gamma}=G_{\gamma}$ or not. 
For $\varphi\in \cc_{c}^{\infty}(\kg)$, $y\in G$, let 
$$
\varphi^{y}(X)=\varphi(\Ad(y)X),\quad \forall\, X\in \kg.  
$$
Following the same lines as \cite{a1}, Lem. 8.1, one can show
\begin{equation}\label{arthur var}
J_{M}(\gamma,\varphi^{y})-J_{M}(\gamma,\varphi)=\sum_{\substack{Q\in \cf(M)\\ Q\neq G}}
J_{M}^{M_{Q}}(\gamma,\varphi_{Q,y}),
\end{equation}
where
$$
\varphi_{Q,y}(Y)=\int_{\kn_{Q}}\int_{K}
\varphi\big(\Ad(k^{-1})(Y+n)\big)\rv_{Q}'(ky)dkdn,\quad \forall\,Y\in \km_{Q}.
$$
Despite the complexity of the expression, the function $\varphi_{Q,y}$ has the nice property:
\begin{equation}\label{key observ}
\text{If }\varphi\in \cc_{c}(\kg/\kk_{r}), \text{ then } \varphi_{Q,y}\in \cc_{c}(\km_{Q}/\km_{Q}\cap \kk_{r}), \quad \forall r\in \br,
\end{equation}
because $\kk_{r}$ is stable under the adjoint action of $K$. 

\subsection{The Shalika germ expansion}

We make a brief recall on the Shalika germ expansion following Kottwitz \cite{kottwitz course}. 
By definition, a \emph{distribution} on $\kg$ is a linear functional on $\cc_{c}^{\infty}(\kg)$. Let $\cd(\kg)$ be the space of distributions on $\kg$, it is endowed with the topology as explained in \cite{kottwitz course}, \S26.3.
Let $\cd(\kg)^{G}$ be the subspace of $G$-invariant distributions on $\kg$. For any  $X\in \kg$, it is known that the orbital integral
$$
J_{G}(X, f)=|D(X)|^{\frac{1}{2}}\int_{G_{X}\backslash G} f(\Ad(g^{-1})X)d\bar{g}
$$
is absolutely convergent for any $f\in \cc_{c}^{\infty}(\kg)$, hence it defines a $G$-invariant distribution $\co_{X}$. Let $\cd(\kg)_{\rm orb}$ be the subspace of $\cd(\kg)^{G}$ spanned by the orbital integral $\co_{X},X\in \kg$. 
It is well known that the subspace $\cd(\kg)_{\rm orb}$ has dense image in $\cd(\kg)^{G}$. 

\begin{thm}[Shalika]

Let $\bg$ be a connected semisimple algebraic group over $F$. There exists functions $\Gamma_{G}(\,\cdot,u)$ on $\kg^{\rm rs}$, $u\in [\nn_{G}(F)]$, satisfying the property:   
For any $f\in \cc_{c}^{\infty}(\kg)$, there exists an open neighbourhood $U_{f}$ of $0$ in $\kg$ such that
$$
\co_{X}(f)=\sum_{u\in [\nn_{G}(F)]}\Gamma_{G}(X,u)\co_{u}(f),\quad \forall\, X\in U_{f}\cap \kg^{\rm rs}. 
$$
The germs around $0$ of the functions $\Gamma_{G}(\,\cdot,u)$ are uniquely determined, called \emph{Shalika germs}.

\end{thm}

The nilpotent orbits in $\kg$ have a nice property: they are preserved under scaling. 
For $u\in \nn_{G}(F)$, let $\co_{u}$ be its $G(F)$-orbit.
It is known that if $t\in F^{\times}$, then $u$ and $tu$ are conjugate over $\overline{F}$; and if $t\in (F^{\times})^{2}$, they will be $G(F)$-conjugate. Moreover, 
\begin{equation}\label{uni meas scale}
d(tx)=|t|^{\dim(\co_{u})/2} dx
\end{equation}
for the canonical measure $dx$ on $\co_{u}$. By the uniqueness of the Shalika germ expansion, we get the homogeneous property of the Shalika germs 
$$
\Gamma_{G}(tX,u)=|t|^{\frac{1}{2}\kparen{\dim(G_{u})-\rk(G)}} \Gamma_{G}(X,u). 
$$
This reduces the study of the orbital integrals to elements with small root valuations, which is much easier. The problem arises then to determine the exact domain of validity of the Shalika germ expansion.

For a compact subset $\omega$ of $\kg$, we denote 
$$
{}^{G}\omega:=\vparen{\Ad(g)X\mid g\in G, X\in \omega}.
$$
For a $G$-invariant subset $\Omega\subset \kg$, we denote by $J(\Omega)$ the invariant distributions on $\kg$ with support in $\Omega$. The following theorem determines the exact domain of validity of the Shalika germ expansion.

\begin{thm}[DeBacker \cite{debacker}]\label{debacker main}

For any $r\in \br$, we have equality of distributions:
$$
\mathrm{Res}_{\ch_{r}}\big(J(\kg_{r})\big)=\mathrm{Res}_{\ch_{r}}\big(J(\nn_{G})\big).
$$

\end{thm}

In general, we can consider the germ expansion around any element $s\in \ka_{M_{\min}}$.

\begin{thm}\label{shalika descent}

Let $\bg$ be a connected reductive algebraic group over $F$, let $s\in \ka_{M_{\min}}$. There exists functions $\Gamma_{G}(s+\cdot,s+u)$ on $\kg_{s}^{\der,\rm rs}$, $u\in [\nn_{G_{s}}(F)]$, satisfying the property:   
For any $f\in \cc_{c}^{\infty}(\kg)$, there exists an open neighbourhood $V_{f}$ of $0$ in $\kg_{s}^{\der}$ such that
$$
J_{G}(s+Y,f)=\sum_{u\in [\nn_{G_{s}}(F)]}\Gamma_{G_{s}}(s+Y,s+u)J_{G}(s+u, f),\quad \forall \,Y\in V_{f}\cap \kg_{s}^{\der,\rm rs}. 
$$

\end{thm}

\begin{proof}

Take a parabolic subgroup $Q\in \cp(G_{s})$, recall that for $Y\in \kg^{\der}_{s},\,f\in \cc_{c}^{\infty}(\kg)$, we have the parabolic descent
\begin{equation}\label{orb int descent}
J_{G}(s+Y,f)=J_{G_{s}}^{G_{s}}(s+Y, f^{Q}), 
\end{equation}
where
\begin{equation}\label{def fq}
f^{Q}(Y')=\int_{K}\int_{\kn_{Q}} f(\Ad(k)(Y'+\nu))d\nu dk, \quad \forall\, Y'\in \kg_{s}. 
\end{equation}
Applying the Shalika germ expansion for $G_{s}$, there exists an open neighbourhood $V_{f}$ of $0$ in $\kg_{s}^{\der}$ such that
\begin{align}
J_{G}(s+Y,f)&=J_{G_{s}}^{G_{s}}(s+Y, f^{Q})=\sum_{u\in [\nn_{G_{s}}(F)]} 
\Gamma_{G_{s}}(s+Y,s+u) J_{G_{s}}^{G_{s}}(s+u, f^{Q})\nonumber
\\
&=\sum_{u\in [\nn_{G_{s}}(F)]} 
\Gamma_{G_{s}}(s+Y,s+u) J_{G}(s+u, f). \label{para descent shalika}
\end{align}

\end{proof}

\begin{thm}\label{debacker descent}

Let $s\in \ka_{M_{\min}}$, we have the equality of distributions
$$
{\rm Res}_{\ch_{r}}J({}^{G}(s+(\kg^{\rm der}_{s})_{r}))={\rm Res}_{\ch_{r}}J({}^{G}(s+\nn_{G_{s}})),\quad\forall\,r\in\br.
$$

\end{thm}

\begin{proof}

According to Harish-Chandra, every invariant distribution on $\kg$ lies in the closure of the linear span of the subset $\vparen{\co_{X}\mid X\in \kg^{\rss}}$. (We refer the read to \cite{kottwitz course}, \S27, for a beautiful proof of this result.) Consequently, $J({}^{G}(s+(\kg_{s}^{\der})_{r}))$ lies in the closure of the linear span of the subset
$$
\vparen{\co_{X}\,\big| \, X\in \kg^{\rss}\cap {}^{G}(s+(\kg_{s}^{\der})_{r})}.
$$

For $Y \in (\kg^{\rm der}_{s})_{r}$, we have parabolic descent as before:
\begin{equation*}
J_{G}(s+Y,f)=J_{G_{s}}^{G_{s}}(s+Y, f^{Q}). 
\end{equation*}
Notice that if $f\in \cc_{c}(\kg/\kk_{r})$, then $f^{Q}\in \cc_{c}(\kg_{s}/\kg_{s}\cap\kk_{r})$ as $\kk_{r}$ is stable under the adjoint action of $K$. 
Hence the expansion (\ref{para descent shalika}) holds for $f\in \cc_{c}(\kg/\kk_{r})$ if we apply theorem \ref{debacker main} to the group $\bg_{s}$. The expansion, together with the above density result, implies the inclusion of distributions
$$
{\rm Res}_{\cc_{c}(\kg/\kk_{r})}J({}^{G}(s+(\kg^{\rm der}_{s})_{r}))\subset {\rm Res}_{\cc_{c}(\kg/\kk_{r})}J({}^{G}(s+\nn_{G_{s}})).
$$
The inclusion in the other direction is clear, hence the equality of distributions as claimed.

\end{proof}

\subsection{The Arthur-Shalika expansion}

For this section, we work with the group $\bg=\gl_{d+1}$. Let $\ba$ be the maximal torus of $\bg$ of the diagonal matrices, let $\mathbf{B}_{0}$ be the Borel subgroup of $\bg$ of the upper triangular matrices. Let $K=\bg(\co)$, let $I$ be the standard Iwahori subgroup, it is the pre-image of $\mathbf{B}_{0}(k)$ under the reduction map $\bg(\co)\to \bg(k)$. In terms of Moy-Prasad filtration, we have $K=G_{0,0}$ and $I=G_{x_{0},0}$ with $x_{0}=\big(1,\frac{d}{d+1},\cdots,\frac{1}{d+1}\big)$. For $n\in \bz$, let $\kk_{n}=\kg_{0,n}$ and $\ki_{n}=\kg_{x_{0},n}$, and similarly for $\kk_{n+}$ and $\ki_{n+}$.

To begin with, we can reinforce DeBacker's homogeneity result. 
Let $\cl_{n+}^{G}$ be the vector space of linear functionals $L:\cc_{c}(\kg)\to \bc$ satisfying the conditions:
\begin{enumerate}[nosep, label=(\roman*)]
\item

For $f\in \cc_{c}(\kg/\kk_{n+})$, we have $L(f^{g})=L(f)$ for any $g\in G$, 

\item

$L(f)=0$ for $f\in \cc_{c}(\kg/\kk_{n+})$ with $\supp(f)\cap (\kk_{n+}+\nn_{G})=\emptyset$. 

\end{enumerate} 

\begin{thm}\label{debacker gln general}

For any $n\in \bz$, we have
$$
\cl_{n+}^{G}=\res_{\cc_{c}(\kg/\kk_{n+})}(J(\nn_{G})).
$$

\end{thm}

\begin{proof}

After scaling by $\ep^{n}$ on $\kg$, the assertion is reduced to the $n=0$ case, which is essentially the main result of DeBacker \cite{debacker gl}. Indeed, the key lemma 4 of \emph{loc. cit.} holds more generally for the linear functionals in $\cl_{0+}^{G}$. We don't need the much stronger assumption that $D$ is an invariant \emph{distribution}, because the test functions involved in the proof lie in $\cc_{c}(\kg/\kk_{0+})$.

\end{proof}

\begin{rem}

The theorem is quite special to the group $\gl_{d+1}$. Although DeBacker \cite{debacker} has established the homogeneity of Shalika expansion for general reductive groups, it doesn't seem that a strengthened result as above will hold in the general setting.   

\end{rem}

\begin{thm} \label{shalika}

Let $M\in \cl(A)$. For each element $u\in \left[\nn_{L}(F)\right], L\in \cl(M)$, there exists a unique function $g_{M}^{L}(\,\cdot\,,u)$ on $\km$, such that for any $ f\in\cc_{c}(\kg/\kk_{n+})$ and any element $\gamma\in  {}^{G}\ki_{n+}\cap \km$, $n\in \bz$, we have
$$
J_{M}(\gamma,f)=\sum_{L\in \cl(M)}\sum_{u\in \left[\nn_{L}(F)\right]}
g_{M}^{L}(\gamma,u)J_{L}(u,f).
$$

\end{thm}

\begin{proof}
Up to scaling, it is enough to prove the assertion for $n=0$. 
We will prove it in this case by induction on $\rk(G)$. For $G=M$, the theorem is exactly the Shalika germ expansion, combined with the theorem of DeBacker.

Now suppose that the theorem holds for the proper Levi subgroups of $G$ containing $M$, then 
\begin{align}
J_{M}(\gamma,f^{y})-J_{M}(\gamma,f)&=\sum_{\substack{Q\in \cf(M)\\Q\neq G}}  J_{M}^{M_{Q}}(\gamma,f_{Q,y}) \nonumber
\\
&=\sum_{\substack{Q\in \cf(M)\\Q\neq G}} \sum_{L\in \cl^{M_{Q}}(M)}\sum_{u\in \kparen{\nn_{L}(F)}} g_{M}^{L}(\gamma,u)\cdot J_{L}^{M_{Q}}(u,f_{Q,y}) \nonumber
\\
&=\sum_{\substack{L\in \cl(M)\\L\neq G}} \sum_{u\in \kparen{\nn_{L}(F)}} g_{M}^{L}(\gamma,u)\sum_{\substack{Q\in \cf(L)\\Q\neq G}}  J_{L}^{M_{Q}}(u,f_{Q,y}) \nonumber
\\
&=\sum_{\substack{L\in \cl(M)\\L\neq G}} \sum_{u\in \kparen{\nn_{L}(F)}} g_{M}^{L}(\gamma,u)\kparen{J_{L}(u,f^{y})-J_{L}(u,f)}, \label{J invariance}
\end{align}
where the first and the last equation is an application of the equation (\ref{arthur var}) to the element $\gamma$ and $u$ respectively, and the second equation is just the induction hypothesis applied to the function $f_{Q,y}$, which belongs to $\cc_{c}(\km_{Q}/\km_{Q}\cap\kk_{0+})$ as we have observed in (\ref{key observ}). 
The above equation implies that the linear functional on $\cc_{c}(\kg/\kk_{0+})$ sending $f\in \cc_{c}(\kg/\kk_{0+})$ to
\begin{equation}\label{def jj 1}
J'_{M}(\gamma,f):=J_{M}(\gamma,f)-\sum_{\substack{L\in \cl(M)\\L\neq G}}\sum_{u\in \left[\nn_{L}(F)\right]}
g_{M}^{L}(\gamma,u)J_{L}(u,f)
\end{equation}
is $G$-conjugation invariant. Moreover, as ${}^{G}\ki_{0+}\subset \kk_{0+}+\nn_{G}$, we have $J_{M}(\gamma,f)=0$ for $f\in \cc_{c}(\kg/\kk_{0+})$ with $\supp(f)\cap (\kk_{0+}+\nn_{G})=\emptyset$, and similarly for $J_{M}(u,f), u\in [\nn_{G}(F)],$ by taking limits. Hence $J_{M}'(\gamma,\,\cdot\,)\in \cl_{0+}^{G}$, and so it is a linear combination of the invariant distributions $\co_{u},u\in [\nn_{G}(F)]$ by theorem \ref{debacker gln general}, this defines the functions $g_{M}^{G}(\gamma,u)$. Their uniqueness follows from the linear independence of the invariant distributions $\co_{u}$'s.

\end{proof}

\begin{rem}

As pointed out by Prof. Waldspurger, the functional $J_{M}'(\gamma, \bullet)$ is not a $G$-invariant \emph{distribution}. Indeed, in the deduction of (\ref{J invariance}), to use the induction hypothesis, $f$ is constrained to lie in $\cc_{c}(\kg/\kk_{0+})$ and it can not be an arbitrary element in $\cc_{c}^{\infty}(\kg)$. 
This is the main obstacle for the generalization of the theorem to arbitrary reductive groups.

\end{rem}

The functions $g_{M}^{G}$ are called the \emph{Arthur-Shalika} germs. In \cite{a1}, Arthur has shown that the \emph{germs} of these functions satisfies a homogeneous property and a semisimple descent formula. 
Now that we have removed all the ambiguities on the Arthur-Shalika germ expansion in \emph{loc. cit.} and determine its domain of validity, we can reinforce these results to the \emph{functions} themselves. For $u\in \nn_{M}(F)$, Arthur has defined a $(G,M)$-family
$$
\rc_{P}(\lambda; \exp(u),t), \quad {P\in \cp(M)}. 
$$
It depends only on the geometric conjugacy class of $u$. By abuse of notation, we will denote it by $\rc_{P}(\lambda; u,t)$. For $v\in [\nn_{L}(F)]$, $L\in \cl(M)$, let 
$$
[u_{M}^{L}:v]=\begin{cases}
1,& \text{ if } v \in u_{M}^{L},\\
0,& \text{ otherwise}, 
\end{cases}
$$
and $d^{M}(u)=\frac{1}{2}(\dim(M_{u})-\rk(M))$. Obviously $d^{M}(u)=d^{L}(v)$ if $v\in u_{M}^{L}$.

\begin{thm}\label{shalika homogeneity}

For $u\in [\nn_{G}(F)]$ and $\gamma\in \km\cap \kg^{\rss}$, we have
$$
g_{M}^{G}(t\gamma,u)=\norm{t}^{d^{G}(u)}\sum_{L\in \cl(M)} \sum_{v\in [\nn_{L}(F)]} g_{M}^{L}(\gamma,v)\rc_{L}(v,t)[v_{L}^{G}:u],\quad \forall \, t\in F^{\times}.
$$

\end{thm}

\begin{proof}

The proof proceeds along the same lines as in \cite{a1}, \S10, taking into account the fact that $g_{M}^{G}(t\gamma,u)=g_{M}^{G}(t\gamma,tu)$ as $u$ is conjugate to $tu$ for $\gl_{d+1}$, and that the scaling on $\kg$ sends
$$
{}^{G}\ki_{n+}\stackrel{{\cdot}^{t}}{\longmapsto} {}^{G}\ki_{n+\val(t)+},\quad \cc_{c}(\kg/\kk_{n+})\stackrel{{\cdot}^{t}}{\longmapsto} \cc_{c}(\kg/\kk_{n+\val(t)+}),\quad \forall n\in \bz,
$$
which preserves the validity of theorem \ref{shalika}.

\end{proof}

\subsection{The parabolic descent}

Let $\bg=\gl_{d+1}$ as before.
Let $\chi:\kg\to \kg\sslash G$ be the characteristic morphism. The above results describe the local behavior of $J_{M}(\gamma,f)$ as $\chi(\gamma)\to 0$. 
We can describe its local behavior as $\chi(\gamma)$ approaches a wall. 
We fix an element $s\in \ka$, and we consider the Levi subgroups $M\in \cl(A)$ such that $s\in \kz_{M}$ or equivalently $M\subset G_{s}$. 
Let $\mu\in \km\cap \kg_{s}^{\der}$ such that $G_{s+\mu}\subset G_{s}$. Following Arthur \cite{a1}, we can express $J_{M}(s+\mu, f)$ in terms of weighted orbital integrals of $\mu$ for the group $G_{s}$. To begin with, we write
\begin{align}
J_{M}(s+\mu,f)&=\lim_{\substack{a\in \ka_{M}^{\reg}\\ a\to 0}} \sum_{L\in \cl(M)} \rr_{M}^{L}(s+\mu, a) J_{L}(a+s+\mu,f)\nonumber
\\
&=\lim_{\substack{a\in \ka_{M}^{\reg}\\ a\to 0}} 
\norm{D(a+s+\mu)}^{\frac{1}{2}}
\int_{G_{s}\backslash G} dg_{1} \int_{G_{a+s+\mu}\backslash G_{s}} dg_{2}\cdot f\big(\Ad(g_{2}g_{1})^{-1} (a+s+\mu)\big)   \nonumber
\\
&\qquad\cdot \kparen{\sum_{L\in \cl(M)} \rr_{M}^{L}(s+\mu, a) \rv_{L}(g_{2}g_{1})}. \label{levi descent 1}
\end{align}

\begin{lem}[Arthur \cite{a1}]
For $L\in \cl(M)$, we have
$$
\rr_{M}^{L}(s+\mu,a)=\begin{cases}
\rr_{M}^{L}(\mu,a),&\text{ if } L\subset G_{s},\\
0,& \text{ otherwise}.
\end{cases}
$$

\end{lem}

For the weight factor $\rv_{L}(g_{2}g_{1})$, it can be written as the volume of the product of two $(G,L)$-families. For $Q=LN_{Q}\in \cp(L)$, we can write $g_{2}=lnk_{Q}(g_{2})$ with $l\in L(F)$, $n\in N_{Q}(F)$, $k_{Q}(g_{2})\in K$, then
\begin{align*}
\rv_{Q}(\lambda; g_{2}g_{1})&=e^{\pair{\lambda,\, H_{Q}(g_{2}g_{1})}}=e^{\pair{\lambda,\, H_{Q}(g_{2})}} e^{\pair{\lambda,\,H_{Q}(k_{Q}(g_{2})g_{1})}}
\\
&=\rv_{Q}(\lambda; g_{2})\ru_{Q}(\lambda;g_{2},g_{1}),
\end{align*}
with $\ru_{Q}(\lambda;g_{2},g_{1})=e^{\pair{\lambda,\,H_{Q}(k_{Q}(g_{2})g_{1})}}$. Notice that  $\big(\ru_{Q}(\lambda;g_{2},g_{1})\big)_{Q\in \cp(L)}$ forms a $(G,L)$-family, hence by the product formula (\ref{arthur prod}), we have
\begin{equation}\label{vol prod}
\rv_{L}(g_{2}g_{1})=\sum_{Q'\in \cf(L)} \rv_{L}^{Q'}(g_{2})\ru_{Q'}'(g_{2},g_{1}).
\end{equation}
By construction, it is clear that for any $k\in K$, we have
$$
\ru_{Q}(\lambda;g_{2},kg_{1})=\ru_{Q}(\lambda;g_{2}k,g_{1}), 
$$
hence 
\begin{equation}\label{uq prime k move}
\ru_{Q'}'(g_{2},kg_{1})=\ru_{Q'}'(g_{2}k,g_{1}),\quad \forall\, Q'\in \cf(L).
\end{equation}

We can further simplify the formula with the condition $g_{2}\in G_{s}$. By the above formula, it is enough to consider the Levi subgroups $L\in \cl^{G_{s}}(M)$. 

\begin{lem}[Arthur \cite{a1}]

For $L\in \cl^{G_{s}}(M)$, we have
$$
\rv_{L}(g_{2}g_{1})=\sum_{R\in \cf^{G_{s}}(L)} \rv_{L}^{R}(g_{2})\ru_{R}'(g_{2},g_{1}),
$$
where
$$
\ru_{R}'(g_{2},g_{1})=\sum_{\{Q'\in \cf(M)\mid Q'_{s}=R,\, \ka_{Q'}=\ka_{R}\}} 
\ru_{Q'}'(g_{2},g_{1}).
$$

\end{lem}

Combining the above two lemmas, we get
$$
\sum_{L\in \cl(M)} \rr_{M}^{L}(s+\mu, a) \rv_{L}(g_{2}g_{1})
=\sum_{R\in \cf^{G_{s}}(M)} \sum_{J\in \cl^{M_{R}}(M)} \rr_{M}^{J}(\mu,a)  \rv_{J}^{R}(g_{2})\ru_{R}'(g_{2},g_{1}).
$$
Plug it into the equation (\ref{levi descent 1}), we get
\begin{align}
&J_{M}(s+\mu,f)
=\lim_{\substack{a\in \ka_{M}^{\reg}(F)\\ a\to 0}} 
\norm{D(a+s+\mu)}^{\frac{1}{2}}
\int_{G_{s}\backslash G} dg_{1} \int_{G_{a+s+\mu}\backslash G_{s}} dg_{2}\cdot f\big(\Ad(g_{2}g_{1})^{-1} (a+s+\mu)\big)\nonumber
\\
&\hspace{3cm}\cdot \kparen{\sum_{R\in \cf^{G_{s}}(M)} \sum_{J\in \cl^{M_{R}}(M)} \rr_{M}^{J}(\mu,a)  \rv_{J}^{R}(g_{2})\ru_{R}'(g_{2},g_{1})} \nonumber
\\
&=\lim_{\substack{a\in \ka_{M}^{\reg}(F)\\ a\to 0}} 
\norm{D(a+s+\mu)}^{\frac{1}{2}}
\int_{G_{s}\backslash G} dg_{1}
\sum_{R\in \cf^{G_{s}}(M)} 
\int_{M_{s+\mu}\backslash M_{R}} dm \int_{N_{R}} dn \int_{K^{G_{s}}} dk \nonumber
\\
&\qquad\cdot f\big(\Ad(g_{1})^{-1}[s+\Ad(mnk)^{-1}(a+\mu)]\big)
\kparen{\sum_{J\in \cl^{M_{R}}(M)} \rr_{M}^{J}(\mu,a)  \rv_{J}^{R}(m)\ru_{R}'(k,g_{1})}  \nonumber
\\
&=\lim_{\substack{a\in \ka_{M}^{\reg}(F)\\ a\to 0}} 
\norm{D(a+s+\mu)}^{\frac{1}{2}}
\int_{G_{s}\backslash G} dg_{1}
\sum_{R\in \cf^{G_{s}}(M)} \sum_{J\in \cl^{M_{R}}(M)} \rr_{M}^{J}(\mu,a)\int_{M_{\mu}\backslash M_{R}} \nonumber
\\
&\qquad\cdot \kparen{\int_{N_{R}} \int_{K^{G_{s}}} 
f\big(\Ad(g_{1})^{-1}[s+\Ad(mnk)^{-1}(a+\mu)]\big)
  \ru_{R}'(k,g_{1})dkdn } \rv_{J}^{R}(m) dm.  \label{weighted para descent}
\end{align}
For $Y\in \km_{R}$, let
\begin{equation}\label{fr defn}
\widetilde{f}_{R,\,g_{1}}(Y)=\int_{N_{R}} \int_{K^{G_{s}}} 
f\big(\Ad(g_{1})^{-1}[s+\Ad(nk)^{-1}Y]\big)
  \ru_{R}'(k,g_{1})dkdn.
\end{equation}
The above equation can be simplified to
\begin{align}
J_{M}(s+\mu,f)&=\lim_{\substack{a\in \ka_{M}^{\reg}\\ a\to 0}} 
\Big|\frac{D(a+s+\mu)}{D^{G_{s}}(a+\mu)}\Big|^{\frac{1}{2}}
\int_{G_{s}\backslash G} dg_{1}
\sum_{R\in \cf^{G_{s}}(M)} \sum_{J\in \cl^{M_{R}}(M)} \rr_{M}^{J}(\mu,a) J_{J}^{M_{R}}(a+\mu, \widetilde{f}_{R,\,g_{1}}) \nonumber
\\
&=\Big|\frac{D(s+\mu)}{D^{G_{s}}(\mu)}\Big|^{\frac{1}{2}}
\int_{G_{s}\backslash G} dg_{1}
\sum_{R\in \cf^{G_{s}}(M)} 
\lim_{\substack{a\in \ka_{M}^{\reg}\\ a\to 0}}  \sum_{J\in \cl^{M_{R}}(M)} \rr_{M}^{J}(\mu,a) J_{J}^{M_{R}}(a+\mu, \widetilde{f}_{R,\,g_{1}}) \nonumber
\\
&=\Big|\frac{D(s+\mu)}{D^{G_{s}}(\mu)}\Big|^{\frac{1}{2}}
\int_{G_{s}\backslash G} dg_{1}
\sum_{R\in \cf^{G_{s}}(M)} J_{M}^{M_{R}}(\mu, \widetilde{f}_{R,\,g_{1}}), \label{weighted para descent 2}
\end{align}
where in the second line we are able to interchange the order of taking limit and integration because the integration is essentially finite as $\widetilde{f}_{R,\,g_{1}}$ vanishes identically for $g_{1}$ outside of a compact region in $ G_{s}\backslash G$. We summarize the above calculations as

\begin{thm}[Parabolic descent of weighted orbital integrals] \label{semisimple descent}

For $s\in \ka$, $M\in \cl^{G_{s}}(A)$, let $\mu\in \km\cap \kg_{s}^{\der}$ be such that $G_{s+\mu}\subset G_{s}$, we have
$$
J_{M}(s+\mu,f)=\Big|\frac{D(s+\mu)}{D^{G_{s}}(\mu)}\Big|^{\frac{1}{2}}
\int_{G_{s}\backslash G} 
\sum_{R\in \cf^{G_{s}}(M)} J_{M}^{R}(\mu, \widetilde{f}_{R,\,g}) dg.
$$

\end{thm}


Let $U_{s}^{+}\subset \kg_{s}$ be the invariant neighborhood of $0$ consisting of the elements $u$ such that for each eigenspace of $\ad(s)$ over $\kg/\kg_{s}$, the eigenvalues of $\ad(s+u)\ad(s)^{-1}$ have valuation $0$.

\begin{thm}[Parabolic descent of Arthur-Shalika expansion]\label{germ descent}

Let $s\in \ka$ be such that the eigenvalues of $\ad(s)$ on $\kg/\kg_{s}$ all have valuation $r$, let $M\in \cl^{G_{s}}(A)$.
Let $\eta\in \km\cap U_{s}^{+}$ be such that $s+\eta\in \kg^{\rss}$ and $G_{s+\eta}\subset G_{s}$, then for any $f\in \cc_{c}(\kg/\kk_{r})$, we have expansion
$$J_{M}(s+\eta,f)=\sum_{L\in \cl^{G_{s}}(M)}  \sum_{u\in [\nn_{L}(F)]} g_{M}^{L}(\eta,u) 
 J_{L}(s+u,f).
$$

\end{thm}

\begin{proof}

By the homogeneity of the Arthur-Shalika expansion, it is enough to prove the assertion for $r=0$. To begin with, we restrict to the case that $G_{s}$ is a maximal proper Levi subgroup of $G$. We fix a parabolic subgroup $Q\in \cp(G_{s})$.

According to the construction in the proof of theorem \ref{shalika}, the functions $g_{M}^{G}$ are defined inductively as coefficients of the expansion
$$
J_{M}'(s+\eta, f)=\sum_{u\in [\nn_{G}(F)]} g_{M}^{G}(s+\eta,u) J_{G}(u,f), \quad \text{for } f\in \cc_{c}(\kg/\kk). 
$$
To determine them, it is then enough to examine the local behavior of $J_{M}(s+\eta, f)$ for a set of test functions $f\in \cc_{c}(\kg/\kk)$ which forms a dual basis of $J_{G}(u,\,\cdot\,)$, $u\in [\nn_{G}(F)]$.  
Since we are working with $\bg=\gl_{d+1}$, the functions $f$ can be taken to be the characteristic functions $[{}^{\bg(k)}(\ep^{-1}u)+\kk], u\in \nn_{G}(k)$. 
For such $f$ and $\eta\in \km\cap U_{s}^{+}$ such that $s+\eta\in \kg^{\rss}$ and $G_{s+\eta}\subset G_{s}$, by theorem \ref{semisimple descent}, we have
\begin{align}
J_{M}(s+\eta,f)&=\Big|\frac{D(s+\eta)}{D^{G_{s}}(\eta)}\Big|^{\frac{1}{2}}
\int_{G_{s}\backslash G} 
\sum_{R\in \cf^{G_{s}}(M)} J_{M}^{R}(\eta, \widetilde{f}_{R,\,g}) dg \nonumber
\\
&=\delta_{Q}(s)\int_{N_{Q}} 
\sum_{R\in \cf^{G_{s}}(M)} J_{M}^{R}(\eta, \widetilde{f}_{R,\,g}) dg,
\label{wo p descent}
\end{align}
where in the second line we have used the Iwasawa decomposition and the fact that $f$ is invariant under $K$-conjugation, and we have set
$
\delta_{Q}(s)=|\det(\ad(s); \kn_{Q})|.$

Look at the integral $\int_{N_{Q}}  J_{M}^{R}(\eta, \widetilde{f}_{R,\,g}) dg$. With the expression (\ref{fr defn}) and the fact that $f$ is $K$-conjugation invariant, we can manipulate
\begin{align}
&\int_{N_{Q}}dg\int_{M_{s+\eta}\backslash M_{R}}\kparen{\int_{N_{R}} \int_{K^{G_{s}}} 
f\big(\Ad(g)^{-1}[s+\Ad(mnk)^{-1}\eta]\big)
  \ru_{R}'(k,g)dkdn } \rv_{M}^{R}(m) dm \nonumber
\\
 &= \int_{N_{Q}}dg \int_{M_{s+\eta}\backslash M_{R}}\kparen{\int_{N_{R}}  \int_{K^{G_{s}}} 
 f\big(\Ad(kgk^{-1})^{-1}[s+\Ad(mn)^{-1}\eta]\big)
 \ru_{R}'(k,g)dkdn } \rv_{M}^{R}(m) dm  \nonumber
\\
 &= \int_{N_{Q}}dg' \int_{M_{s+\eta}\backslash M_{R}}\kparen{\int_{N_{R}}  \int_{K^{G_{s}}} 
 f\big(\Ad(g')^{-1}[s+\Ad(mn)^{-1}\eta]\big)
 \ru_{R}'(k,k^{-1}g'k)dkdn } \rv_{M}^{R}(m) dm  \nonumber
\\
 &= \int_{N_{Q}}dg' \int_{M_{s+\eta}\backslash M_{R}}\kparen{\int_{N_{R}}   
 f\big(\Ad(g')^{-1}[s+\Ad(mn)^{-1}\eta]\big)
 \ru_{R}'(1,g')dn } \rv_{M}^{R}(m) dm,  \nonumber
\\
 &= \int_{M_{s+\eta}\backslash M_{R}}\int_{N_{R}} \kparen{  
\int_{N_{Q}} f\big(\Ad(g')^{-1}[s+\Ad(mn)^{-1}\eta]\big)
 \ru_{R}'(1,g')dg'  } \rv_{M}^{R}(m) dndm, \label{int jr} 
\end{align}
where the second to the last equality is due to equation (\ref{uq prime k move}) and the right-$K$-invariance of Arthur's weight factor.
The factor $\ru_{R}'(1,g')$ can be further simplified. By construction,
$$
\ru_{R}'(1,g')=\sum_{\{Q'\in \cf(M)\mid Q'_{s}=R,\, \ka_{Q'}=\ka_{R}\}} 
\ru_{Q'}'(1,g')=\sum_{\{Q'\in \cf(M)\mid Q'_{s}=R,\, \ka_{Q'}=\ka_{R}\}} 
\rv_{Q'}'(g').
$$
In particular, it is left-$R$-invariant and right-$K$-invariant in the variable $g'$.

For the integral to be non-zero (and hence positive), we need to solve the equation
\begin{equation}\label{solve springer}
\Ad(g')^{-1}[s+\Ad(mn)^{-1}\eta]\in {}^{\bg(k)}(\ep^{-1}u)+\kk,
\end{equation}  
with $g'\in N_{Q}, m\in M_{R},n\in N_{R}$.
Let $g'=\exp(Z)$ with $Z\in \kn_{Q}$, then the equation can be rewritten
$$
(s+\Ad(mn)^{-1}\eta)+\big[s+\Ad(mn)^{-1}\eta, Z\big]\in {}^{\bg(k)}(\ep^{-1}u)+\kk.
$$
Notice that the left hand side is an element of $\kq$. To solve the equation, we need to understand the intersection ${}^{\bg(k)}u\cap \kq(k)$. It is clear that it is stable under $\bg_{s}(k)$-conjugation, hence there exists finitely many orbits ${}^{\bg_{s}(k)}u_{i}$ in $\nn_{G_{s}}(k)$, ${i\in I}$, such that 
$$
{}^{\bg(k)}u\cap \kq(k)=\bigsqcup_{i\in I}  {}^{\bg(k)}u  \cap \paren{{}^{\bg_{s}(k)}u_{i}+\kn_{Q}}.
$$
Hence the set of solutions $\Omega$ to the equation (\ref{solve springer}) is the disjoint union of the set of solutions $\Omega_{i}$ to the equations
\begin{align}
&\Ad(mn)^{-1}\eta\in {}^{\bg_{s}(k)}(\ep^{-1}u_{i})+\kk^{G_{s}},  \label{sol levi}
\\
&\big[s+\Ad(mn)^{-1}\eta, Z\big]\in \ep^{-1}(\kn_{Q}\cap \kk),  \label{solve nilp}
\end{align}
such that their sum satisfies
\begin{equation} \label{solve total}
\Ad(mn)^{-1}\eta+\big[s+\Ad(mn)^{-1}\eta, Z\big]\in {}^{\bg(k)}(\ep^{-1}u)+\kk. 
\end{equation}

To solve the equation (\ref{sol levi}), note that its left hand side is an element of $\mathfrak{r}$. As before, we can consider the $R(k)$-conjugation and make the decomposition
$$
{}^{\bg_{s}(k)}u_{i}\cap \mathfrak{r}(k)=\bigsqcup_{j\in J_{i}} {}^{R(k)}u_{ij},\quad \text{ with nilpotent } u_{ij}\in \mathfrak{r}(k). 
$$  
Let $\Omega_{ij}$ be the set of solutions to the equation
\begin{equation}\label{solve levi var}
\Ad(mn)^{-1}\eta\in {}^{R(k)}(\ep^{-1}u_{ij})+(\mathfrak{r}\cap \kk), 
\end{equation}
and the equations (\ref{solve nilp}) and (\ref{solve total}), then $\Omega_{i}=\bigsqcup_{j\in J_{i}}\Omega_{ij}$.

Notice that if we conjugate the equation (\ref{solve levi var}) on both sides by an element $h\in R(k)$, we get a solution to the equations (\ref{solve nilp}) and (\ref{solve total}) if we substitute $Z$ by $\Ad(h)Z$.  
Hence we can normalize the equation (\ref{solve levi var}) as 
\begin{equation} 
\label{sol levi var}
\Ad(mn)^{-1}\eta\in \ep^{-1}u_{ij}+\kr\cap \kk.
\end{equation} 
We fix a solution $(m_{0},n_{0},Z_{0})$ of the equation \textup{(\ref{sol levi var})}, \textup{(\ref{solve nilp})} and \textup{(\ref{solve total})}, let $\Ad(m_{0}n_{0})^{-1}\eta=\ep^{-1}u_{ij}+r_{0}$ for some $r_{0}\in \kr\cap \kk$.


\begin{lem} 

For $\eta\in U_{s}^{+}$, let $(m,n,Z)$ be a solution of the equation \textup{(\ref{sol levi var})}, \textup{(\ref{solve nilp})} and \textup{(\ref{solve total})}, let $\Ad(mn)^{-1}\eta=\ep^{-1}u_{ij}+r$ with $r\in \kr\cap \kk$, then there exists $k'\in R\cap K$ such that 
\begin{equation*}
\ad\big(s+\ep^{-1}u_{ij}+r\big)=\ad(s+\ep^{-1}u_{ij}+r_{0}) \circ \Ad(k') \quad \text{on }\kn_{Q}.
\end{equation*}

\end{lem}

\begin{proof}

Suppose that $\ad(s)$ acts on $\kn_{Q}$ with the scalar $t\in \co_{F}^{\times}$.
To begin with, we can reduce to the case that $t=1$. Indeed, change $u_{ij}$ to $tu_{ij}$ does't change the conjugacy class of $u_{ij}$, then change $\eta$ to $t\eta'$, the equations \textup{(\ref{solve levi var})}, \textup{(\ref{solve nilp})} and \textup{(\ref{solve total})} simplifies to the case that $\ad(s)$ acts on $\kn_{Q}$ as identity, with $\eta'$ in the place of $\eta$. But it is clear that $\eta'$ remains in $ U_{s}^{+}$. The question is thus reduced to finding $k'$ such that
\begin{equation*}
1+\ad(\ep^{-1}u_{ij}+r)=(1+\ad(\ep^{-1}u_{ij}+r_{0}))\circ \Ad(k')\quad \text{on }\kn_{Q}.
\end{equation*}

By assumption, $\eta$ is a semisimple integral element, and the eigenvalues of  $1+\ad(\ep^{-1}u_{ij}+r_{0})=1+\ad(\Ad(mn)^{-1}\eta)$ on $\kn_{Q}$ have valuation $0$. 
Hence there exists a semisimple element $k_{0}\in M_{R}$ with eigenvalues of valuation $0$ and a topologically nilpotent element $u_{0}\in \kr$ such that
$$
1+\ad(\ep^{-1}u_{ij}+r_{0})=(1+\ad(u_{0}))\circ \Ad(k_{0})\quad \text{ on }\kn_{Q}.
$$
Similarly, there exists a semisimple element $k_{1}\in M_{R}$ with eigenvalues of valuation $0$,  and a topologically nilpotent element $u_{1}\in \kr$ such that
$$
1+\ad(\ep^{-1}u_{ij}+r)=(1+\ad(u_{1}))\circ \Ad(k_{1})\quad \text{ on }\kn_{Q}.
$$
As $u_{0}, u_{1}$ are topologically nilpotent, we can find elements $A_{0},A_{1}\in \kr$ such that
$$
1+\ad(u_{0})=\exp(\ad A_{0}), \quad 1+\ad(u_{1})=\exp(\ad A_{1}).
$$ 
Let $A=A_{1}-A_{0}\in \kr$, then
\begin{align*}
1+\ad(\ep^{-1}u_{ij}+r)&=(1+\ad(u_{1}))\circ \Ad(k_{1})=(1+\ad(u_{0}))\circ \exp(\ad A)\circ \Ad(k_{1})
\\
&=(1+\ad(\ep^{-1}u_{ij}+r_{0}))\circ\Ad(k_{0}^{-1})\circ \Ad(\exp(A))\circ \Ad(k_{1}).
\end{align*}
Let $k'=k_{0}^{-1}\exp(A)k_{1}\in R$, it belongs to $K$ because $1+\ad(\ep^{-1}u_{ij}+r)$ and $1+\ad(\ep^{-1}u_{ij}+r_{0})$ acts on $\kn_{Q}$ with the same polar part. 
This finishes the proof of the lemma.

\end{proof}


Consequently, the transformation $Z=\Ad(k')Z'$ sets up a bijection between the set of solutions of the equations
\begin{align*}
\big[s+\Ad(mn)^{-1}\eta, Z\big]\in \ep^{-1}(\kn_{Q}\cap \kk),\\
\big[s+\ep^{-1}u_{ij}+r_{0}, Z'\big]\in \ep^{-1}(\kn_{Q}\cap \kk).
\end{align*}
Observe that $\ru_{R}'(1,g')$ is invariant under $\Ad(k')$, as it is left-$R$-invariant and right-$K$-invariant in the variable $g'$. 
Hence for $(m, n)$ satisfying the equation (\ref{solve levi var}), the inner integral in (\ref{int jr}) simplifies as
\begin{align*}
\int_{N_{Q}} f\big(\Ad(g')^{-1}[s+\Ad(mn)^{-1}\eta]\big)
 \ru_{R}'(1,g')dg' 
&= \phi_{R,ij} \int_{N_{Q}} f\big(\Ad(g')^{-1}[s+\Ad(mn)^{-1}\eta]\big)
 dg'
\\
&=\phi_{R,ij} f^{Q}\big(s+\Ad(mn)^{-1}\eta\big)
\end{align*}
for a constant $\phi_{R,ij}$ depending only on the conjugacy class ${}^{R(k)}u_{ij}$.

Let $\Phi^{R}_{ij}$ be the characteristic function of the set ${}^{R(k)}(\ep^{-1} u_{ij})+(\kr\cap \kk)$.
The expression (\ref{wo p descent}) for $J_{M}(s+\eta, f)$ can be simplified as
{
\begin{align}
&J_{M}(s+\eta,f)
=\delta_{Q}(s) \int_{N_{Q}} 
\sum_{R\in \cf^{G_{s}}(M)} J_{M}^{R}(\eta, \widetilde{f}_{R,\,g}) dg \nonumber
\\
&=\delta_{Q}(s)
\sum_{R\in \cf^{G_{s}}(M)} \sum_{i\in I,\, j\in J_{i}}
\int_{M_{s+\eta}\backslash M_{R}}  \int_{N_{R}} \Bigg[\int_{N_{Q}} \Phi^{R}_{ij}(\Ad(mn)^{-1}\eta)\cdot \nonumber
\\
&\hspace{7cm} f\big(\Ad(g)^{-1}[s+\Ad(mn)^{-1}\eta]\big)
  \ru_{R}'(1,g') dn\Bigg] \rv_{M}^{R}(m) dmdg  \nonumber
 \\
&=\delta_{Q}(s)
\sum_{R\in \cf^{G_{s}}(M)} \sum_{i\in I,\, j\in J_{i}} 
\phi_{R,ij}
\int_{M_{s+\eta}\backslash M_{R}}\int_{N_{R}} 
\Phi^{R}_{ij}(\Ad(mn)^{-1}\eta)\cdot  \nonumber
\\
&\hspace{8cm}f^{Q}\big(s+\Ad(mn)^{-1}\eta\big)
 \rv_{M}^{R}(m) dndm   \nonumber
 \\
&=\delta_{Q}(s)
\sum_{R\in \cf^{G_{s}}(M)} \sum_{i\in I,\, j\in J_{i}}
\phi_{R, ij} 
J_{M}^{R}(s+\eta, \Phi^{R}_{ij}\cdot f^{Q}). \label{weighted para eta}
\end{align}
}

Following the same lines of calculations, for general elements $\mu\in \km\cap U_{s}^{+}$, we can simplify the integral $\int_{N_{Q}}J_{J}^{M_{R}}(a+\mu,\widetilde{f}_{R,g_{1}})dg_{1}$ in the expression (\ref{weighted para descent 2}) for $J_{M}(s+\mu,f)$ as:
\begin{align*}
&\int_{N_{Q}}dg'\int_{M_{s+\eta}\backslash M_{R}}\kparen{\int_{N_{R}}   
 f\big(\Ad(g')^{-1}[s+\Ad(mn)^{-1}(a+\mu)]\big)
 \ru_{R}'(1,g')dn } \rv_{J}^{R}(m) dm\\
&=\sum_{i\in I,\, j\in J_{i}}
\phi_{R,ij}
\int_{M_{s+\eta}\backslash M_{R}}\int_{N_{R}} 
\Phi_{ij}^{R}\big(\Ad(mn)^{-1}(a+\mu)\big)\cdot 
 \\
&\hspace{6cm} f^{Q}\big(s+\Ad(mn)^{-1}(a+\mu)\big)\rv_{J}^{R}(m) dndm \\
&=\sum_{i\in I,\, j\in J_{i}}
\phi_{R, ij} 
J_{J}^{M_{R}}(s+a+\mu, \Phi_{ij}^{R}\cdot f^{Q}).
\end{align*}
Plug it into the equation (\ref{weighted para descent 2}) and take limits, we get
\begin{align}
J_{M}(s+\mu,f)&=\Big|\frac{D(s+\mu)}{D^{G_{s}}(\mu)}\Big|^{\frac{1}{2}}
\sum_{R\in \cf^{G_{s}}(M)} 
\lim_{\substack{a\in \ka_{M}^{\reg}\\ a\to 0}}  \sum_{J\in \cl^{M_{R}}(M)} \rr_{M}^{J}(\mu,a) \int_{G_{s}\backslash G}
J_{J}^{M_{R}}(a+\mu, \widetilde{f}_{R,\,g_{1}}) dg_{1} \nonumber
\\
&=\delta_{Q}(s)
\sum_{R\in \cf^{G_{s}}(M)}\lim_{\substack{a\in \ka_{M}^{\reg}\\ a\to 0}}  \sum_{J\in \cl^{M_{R}}(M)} \rr_{M}^{J}(\mu,a)
\sum_{i\in I,\, j\in J_{i}}
\phi_{R,ij} 
J_{J}^{M_{R}}(s+a+\mu, \Phi^{R}_{ij}\cdot f^{Q})    \nonumber
 \\
&=\delta_{Q}(s)
\sum_{R\in \cf^{G_{s}}(M)} \sum_{i\in I,\, j\in J_{i}}
 \phi_{R,ij} J_{M}^{M_{R}}(s+\mu, \Phi^{R}_{ij} \cdot f^{Q}). \label{weighted para mu}
 \end{align}
Now that $\Phi^{R}_{ij}, f^{Q}\in \cc_{c}(\kg_{s}/\kk^{G_{s}})$ and that $s,\eta\in (\kg_{s})_{0}$, we can apply the Arthur-Shalika germ expansion to $J_{M}^{M_{R}}(s+\eta, \Phi^{R}_{ij} \cdot f^{Q})$ and obtain the searched-for expansion
\begin{align}
&J_{M}(s+\eta,f)=\delta_{Q}(s)
\sum_{R\in \cf^{G_{s}}(M)} \sum_{i\in I,\, j\in J_{i}}
\phi_{R,ij} J_{M}^{M_{R}}(s+\eta,\Phi^{R}_{ij} \cdot f^{Q}) \nonumber
\\
&=\delta_{Q}(s)
\sum_{R\in \cf^{G_{s}}(M)} \sum_{i\in I,\, j\in J_{i}}
\phi_{R,ij} 
\sum_{L\in \cl^{M_{R}}(M)}  \sum_{u\in [\nn_{L}(F)]}
g_{M}^{L}(s+\eta,s+u) J_{L}^{M_{R}}(s+u, \Phi^{R}_{ij} \cdot f^{Q})  \nonumber
\\
&=\delta_{Q}(s)
\sum_{L\in \cl^{G_{s}}(M)}  \sum_{u\in [\nn_{L}(F)]} g_{M}^{L}(\eta,u)
\sum_{R\in \cf^{G_{s}}(L)} \sum_{i\in I,\, j\in J_{i}}
\phi_{R,ij} 
 J_{L}^{M_{R}}(s+u, \Phi^{R}_{ij} \cdot f^{Q})     \nonumber
  \\
&=\sum_{L\in \cl^{G_{s}}(M)}  \sum_{u\in [\nn_{L}(F)]} g_{M}^{L}(\eta,u)  J_{L}(s+u,f),  \label{shalika para descent 1}
 \end{align}
where the first and last line are due to equation (\ref{weighted para eta}) and (\ref{weighted para mu}) respectively. This finishes the proof of the theorem when $G_{s}$ is a maximal proper Levi subgroup of $G$.

In general, we can apply the equation (\ref{weighted para eta}) and (\ref{weighted para mu})  iteratively to a sequence of Levi subgroups to conclude.  
Without loss of generality, we can assume that $s=\diag(a_{1}{\rm Id}_{n_{1}},\cdots, a_{l}{\rm Id}_{n_{l}})$. Let $s_{i}=\diag(a_{1}{\rm Id}_{n_{1}},\cdots,a_{i-1}{\rm Id}_{n_{i-1}}, a_{i}{\rm Id}_{n_{i}},\cdots, a_{i}{\rm Id}_{n_{i}})$, then
$$
G_{s}=G_{s_{l}}\subset G_{s_{l-1}}\subset \cdots \subset G_{s_{1}}=G,
$$ 
and each factor $G_{s_{i+1}}$ is a maximal proper Levi subgroup of $G_{s_{i}}$. Take $Q_{i}\in \cp(G_{s_{i}})$ such that
$$
Q_{s}=Q_{l}\subset Q_{l-1}\subset \cdots \subset Q_{1}=G.
$$ 
The equation (\ref{weighted para eta}) implies that 
$$
J_{M}(s+\eta)=J_{M}(s_{2}+(s-s_{2}+\eta))=\delta_{Q_{2}}(s_{2})
\sum_{R'\in \cf^{G_{s_{2}}}(M)} \sum_{i\in I,\, j\in J_{i}}
\phi_{R', ij} 
J_{M}^{M_{R'}}(s+\eta, \Phi^{R'}_{ij}\cdot f^{Q_{2}}).
$$
For $R'\in \cf^{G_{s_{2}}}(M)$ with $M_{R'}\nsubseteq G_{s}$, we can find $i\in \{3,\cdots, l\}$ such that $(M_{R'})_{s_{i}}$ is a maximal proper Levi subgroup of $M_{R'}$, 
and apply parabolic decent (\ref{weighted para eta}) to 
$$
J_{M}^{M_{R'}}(s+\eta, \Phi^{R'}_{ij}\cdot f^{Q_{2}})=J_{M}^{M_{R'}}(s_{i}+(s-s_{i}+\eta), \Phi^{R'}_{ij}\cdot f^{Q_{2}}).
$$ 
This process can be iterated. At the end, we get an expansion of the form
\begin{equation}\label{weighted para eta general}
J_{M}(s+\eta)=\sum_{R\in \cf^{G_{s}}(M)}\sum_{\alpha\in \Upsilon} \varphi_{\alpha}(s)J_{M}^{M_{R}}(s+\eta, f_{\alpha}),
\end{equation}
for some functions $\varphi_{\alpha}(s)$ and $f_{\alpha}\in \cc_{c}(\kg_{s}/\kg_{s}\cap\kk_{r})$, both depending only on $R$ and a nilpotent conjugacy class in $R(k)$, here $\Upsilon$ is some index set. In exactly the same way, we apply the equation (\ref{weighted para mu})  iteratively to $J_{M}(s+u), u\in \nn_{M}(F)$, and get 
\begin{equation}\label{weighted para u general}
J_{M}(s+u)=\sum_{R\in \cf^{G_{s}}(M)}\sum_{\alpha\in \Upsilon} \varphi_{\alpha}(s)J_{M}^{M_{R}}(s+u, f_{\alpha}).
\end{equation}
The looking-for expansion then follows the same line of reasoning as in the deduction of equation (\ref{shalika para descent 1}) from the equations (\ref{weighted para eta general}) and (\ref{weighted para u general}) and the Arthur-Shalika expansion for $J_{M}^{M_{R}}(s+\eta, f_{\alpha})$.

\end{proof}

\section{Counting points on the fundamental domain}\label{counting fund domain}

For this section, we work with  a connected reductive algebraic group $\bfg$ over $\fq$, and for technical simplicity we assume that $\bg^{\der}$ is simply connected. In \cite{chen2}, we have explained a transition between counting points on the fundamental domain of the affine Springer fibers and Arthur's weighted orbital integrals. The strategy is to count the number of points on the truncated affine Springer fibers in two ways: by the Arthur-Kottwitz reduction and by the Harder-Narasimhan reduction. Here we explain another way to count the number of points on the fundamental domains.

Let $T$ be a maximal torus of $G$, let $\gamma\in \kt$ be a regular integral element. Recall that the affine Springer fiber is the sub ind-$k$-scheme of the affine grassmannian $\xx=\bg(\!(\ep)\!)/\bg[\![\ep]\!]$ defined by the equation
$$
\xx_{\gamma}=\big\{g\in \xx \mid \Ad(g)^{-1}\gamma \in \kg[\![\ep]\!]\big\}. 
$$
The group $T$ acts on it by left translation. Let $S\subset T$ be the maximal $F$-split subtorus of $T$, let $\Lambda=\vparen{\ep^{\nu}\mid \nu\in X_{*}(S)}$. Then $\Lambda$ acts freely on $\xx_{\gamma}$ and the quotient $\Lambda\backslash \xx_{\gamma}$ is a projective  $k$-scheme \cite{kl}. Moreover, there is an open dense subscheme $\xx_{\gamma}^{\reg}$ consisting of the points $[g]$ such that $\Ad(g)^{-1}\gamma\mod \ep$ is regular in $\kg$, and the group $T$ acts transitively on $\xx_{\gamma}^{\reg}$ \cite{be}. Let $M_{0}=Z_{\bg_{F}}(S)$ be the centralizer of $S$ in $\bg_{F}$, then $M_{0}$ is a Levi subgroup of $G$ and $T$ is elliptic in $M_{0}$. 
Goresky, Kottwitz and MacPherson \cite{gkm3} have given a characterisation of the regular points in $\xx_{\gamma}$. To formulate it, we need to define an invariant $n(\gamma,P_{1},P_{2})\in \bz_{\ge 0}$ for any two adjacent parabolic subgroups $P_{1}=M_{0}N_{1},P_{2}=M_{0}N_{2}\in \cp(M_{0})$. The Galois group $\gal(\overline{F}/F)$ acts on the set of roots of $T_{\overline{F}}$ in $\kn_{1}\cap \kn_{2}^{-}$. Let $\alpha$ be such a root, let $F_{\alpha}$ be the field of definition of $\alpha$. Let $\val_{F_{\alpha}}$ be the valuation normalised such that any uniformiser in $F_{\alpha}$ has valuation $1$, i.e. $\val_{F_{\alpha}}(\ep)=[F_{\alpha}:F]$. Let $m_{\alpha}$ be the unique positive integer such that the image of $\alpha^{\vee}$ in $\Lambda_{M_{0}}$ is equal to $m_{\alpha}\cdot \beta_{P_{1},P_{2}}$. Now we define
$$
n(\gamma, P_{1},P_{2})=\sum \val_{F_{\alpha}}(\alpha(\gamma))\cdot m_{\alpha},
$$
where the sum is taken over a set of representatives $\alpha$ of the orbits of $\gal(\overline{F}/F)$ on the set of roots of  $T_{\overline{F}}$ in $\kn_{1}\cap \kn_{2}^{-}$.
Recall that for any point $x\in \xx$ the parabolic reductions
$$
\ec_{M_{0}}(x)=\big(H_{P}(x)\big)_{P\in \cp(M_{0})}
$$ 
form a $(G,M_{0})$-orthogonal set, i.e. for any two adjacent parabolic subgroups $P_{1}, P_{2}\in \cp(M_{0})$, there exists $n(x,P_{1},P_{2})\in \bz_{\ge 0}$ such that 
$$
H_{P_{1}}(x)-H_{P_{2}}(x)=n(x,P_{1},P_{2})\cdot \beta_{P_{1}, P_{2}}^{\vee}.
$$

\begin{prop}[Goresky-Kottwitz-MacPherson]\label{gkmreg}

Let $x\in \xx_{\gamma}$. 

\begin{enumerate}[topsep=0pt, itemsep=0pt, label=$(\arabic*)$]

\item For any two adjacent parabolic subgroups $P_{1}, P_{2}\in \cp(M_{0})$, we have
$$
n(x,P_{1},P_{2})\leq n(\gamma,P_{1},P_{2}).
$$

\item The point $x$ is regular in $\xx_{\gamma}$ if and only if the following two conditions holds:
\begin{enumerate}[topsep=0pt, itemsep=0pt, label=\textup{(\roman*)}]
\item
the point $f_{P}(x)$ is regular in $\xx^{M_{0}}_{\gamma}$ for all $P\in \cp(M_{0})$;

\item
for any two adjacent parabolic subgroups $P_{1}, P_{2}$ in $\cp(M_{0})$, one has
$$
n(x,P_{1},P_{2})=n(\gamma,P_{1},P_{2}).
$$

\end{enumerate}
\end{enumerate}
\end{prop}

This beautiful result motivates us to introduce the \textit{fundamental domain} $F_{\gamma}$ of $\xx_{\gamma}$:  
Fix a point $x_{0}\in \xx_{\gamma}^{\reg}$, we define
$$
F_{\gamma}=\Big\{x\in \xx_{\gamma}\,\Big|\, \ec_{M_{0}}(x)\subset \ec_{M_{0}}(x_{0}),\, x  \text{ and } x_{0} \text{ lie on the same connected component}\Big\}.
$$
Then $F_{\gamma}$ is a projective $k$-scheme, and for different choices of $x_{0}$, the fundamental domains are isomorphic. As explained in \cite{chen1}, the study of the geometry of $\xx_{\gamma}$ can be reduced to that of $F_{\gamma}$ by the Arthur-Kottwitz reduction.

From now on we assume that the splitting field of $T$ is totally ramified over $F$. Let $T^{1}=T\cap \ker (H_{M_{0}})$. By \cite{chen2}, lemma 2.6, $T^{1}$ is of finite volume and we have an exact sequence
\begin{equation}\label{T devi}
1\to T^{1}\to T\xrightarrow{H_{M_{0}}} X_{*}(M_{0})\to 1.
\end{equation}
Let $\Lambda_{M_{0}}$ be the quotient of $X_{*}(A)$ by the coroot lattice of $M_{0}$, recall that the connected components of $\xx_{\gamma}^{M_{0}}$ are parametrized by $\Lambda_{M_{0}}$. Now that $M_{0}^{\rm der}$ is also simply connected, we have $X_{*}(M_{0})=\Lambda_{M_{0}}$ by \cite{cl2}, lemma 11.6.1. For any point $x\in \xx_{\gamma}$, we define a weight factor
$$
\rv_{\gamma}(x)=|\vparen{\lambda\in \Lambda_{M_{0}}\mid \lambda+\ec_{M_{0}}(x)\subset \ec_{M_{0}}(x_{0})}|. 
$$

\begin{lem}\label{fund weight orb}

We have 
$$
|F_{\gamma}(\fq)|=\vol(T^{1})\int_{T\backslash G} \mathbbm{1}_{\kk}(\Ad(g^{-1})\gamma)\rv_{\gamma}(g)d\bar{g}.
$$

\end{lem}

\begin{proof}

Take a parabolic subgroup $P=M_{0}N\in \cp(M_{0})$, consider the parabolic reduction
$$
f_{P}:\xx_{\gamma}\to \xx_{\gamma}^{M_{0}}.
$$ 
It is known that its restriction to the pre-image of each connected component $\xx_{\gamma}^{M_{0},\nu}$, $\nu\in \Lambda_{M_{0}}$, is an iterated fibration in affine spaces, and that the translation with $T$ identifies them \cite{kl}. In particular, every point in $F_{\gamma}$ can be moved to $f_{P}^{-1}\big(\xx_{\gamma}^{M_{0},0}\big)$ by translation with $T$. Combined with the fact that $H_{P}=H_{M_{0}}\circ f_{P}$, and that $T$ fits into the exact sequence (\ref{T devi}), we get
\begin{align*}
|F_{\gamma}(\fq)|&=\sum_{x\in f_{P}^{-1}\big(\xx_{\gamma}^{M_{0},0}\big)(\fq)}
|\vparen{\lambda\in \Lambda_{M_{0}}\mid \lambda+\ec(x)\subset \ec(x_{0})}|
\\
&=\vol(T^{1})\int_{T\backslash G} \mathbbm{1}_{\kk}(\Ad(g^{-1})\gamma)\rv_{\gamma}(g)d\bar{g}.
\end{align*}

\end{proof}

For $x\in \xx_{\gamma}$, let $\lambda(x,x_{0})$ be the positive $(G,M_{0})$-orthogonal set
$$
\lambda_{P}(x,x_{0})=H_{P}(x_{0})-H_{P}(x),\quad \forall\,P\in \cp(M_{0}).
$$
In other words, it is the difference of the $(G,M_{0})$-orthogonal sets $\ec_{M_{0}}(x_{0})$ and $\ec_{M_{0}}(x)$. Here the positivity of $\lambda(x,x_{0})$ follows from the characterization of the regular points on $\xx_{\gamma}$ by Goresky, Kottwitz and MacPherson \cite{gkm3}.

\begin{lem}\label{fund weight factor}

For any point $x\in \xx_{\gamma}$, we have
$$
\rv_{\gamma}(x)=|\be(\lambda(x,x_{0}))\cap \Lambda_{M_{0}}|.
$$

\end{lem}

\begin{proof}

By definition, the inclusion of polytopes $\lambda+\ec(x)\subset \ec(x_{0})$ is equivalent to the condition
$$
\lambda+H_{P}(x)\prec_{P} H_{P}(x_{0}),\quad \forall\,P\in \cp(M_{0}),
$$
where $a\prec_{P_{0}}b$ means that $b-a$ is a positive linear combination of elements in $\Delta_{P}^{\vee}$. Obviously, this is equivalent to the condition that $\lambda\in \be(\lambda(x,x_{0}))$.

\end{proof}

To count the number of lattice points in the polytope generated by a positive $(G,M)$-orthogonal set, the question has been addressed in \cite{chen2}, \S2.2.3, following the work   of Arthur \cite{a4} and Chaudouard-Laumon \cite{cl2}. Let $h=(h_{P})_{P\in \cp(M)}$ be an \emph{integral} positive $(G,M)$-orthogonal set. By the word \emph{integral}, we mean that $h_{P}\in \Lambda_{M}$ for all $P\in \cp(M)$.
Let $\rw_{M}(h)$ be the number of lattice points in $\Lambda_{M}$ contained in the polytope $\be(h)$, let $\rv_{M}(h)$ be the volume of the polytope $\be(h)$.

\begin{lem}\label{arthur lattice}
Let 
$
\overline{\rc}_{P}(\lambda)=\prod_{\alpha\in \Delta_{P}} (1-e^{-\langle\lambda,\,\alpha^{\vee}\rangle}),
$
then
\begin{equation*}
\rw_{M}(h)=
\lim_{\lambda\to 0}
\sum_{P\in \cp(M)}
\overline{\rc}_{P}(\lambda)^{-1}e^{\langle \lambda,\, h_{P}  \rangle}.
\end{equation*}

\end{lem}

\begin{proof}

The proof follows the same line as the deduction of the equations (2.8) in \cite{chen2}, we indicates only the modifications. Let ${\varphi}_{P}^{\lambda}$ be the characteristic function of the cone
$$
\{a\in \ka_{M}\mid \varpi_{\alpha}(a)> 0,\,\forall \,\alpha\in \Delta_{P}^{\lambda};\;\varpi_{\alpha}(a)\le 0,\,\forall \,\alpha\in \Delta_{P}\backslash \Delta_{P}^{\lambda}\},
$$
where $\lambda$ is a generic point in $\ka_{M}^{*}$, and 
$$
\Delta_{P}^{\lambda}=\{\alpha\in \Delta_{P}\mid \langle\lambda, \alpha\rangle<0\}.
$$ 
Then the characteristic function of the \emph{closed} convex polytope $\be(h)$ is equal to the function
$$
a\in \ka_{M}\longmapsto \sum_{P\in \cp(M)} (-1)^{|\Delta_{P}^{\lambda}|}{\varphi}_{P}^{\lambda}(-h_{P}+a),
$$
and the function $\rw_{M}(h)$ can be expressed as
$$
\rw_{M}(h)=
\sum_{\chi\in X_{*}(M)}\sum_{P\in \cp(M)}(-1)^{|\Delta_{P}^{\lambda}|} {\varphi}^{\lambda}_{P}(-h_{P}+\chi).
$$
Again we introduce an extra exponential factor to treat the infinite sum in $X_{*}(M)$, let 
$$
\overline{S}_{P}(\lambda)=\sum_{\chi\in X_{*}(M)} {\varphi}^{\lambda}_{P}(-h_{P}+\chi)e^{\langle\lambda,\chi\rangle}.
$$
The series converges absolutely for generic $\lambda$, and
$$
\rw_{M}(h)=\lim_{\lambda\to 0}\sum_{P\in \cp(M)}(-1)^{|\Delta_{P}^{\lambda}|} \overline{S}_{P}(\lambda),
$$
where the limit is taken for generic $\lambda\in \ka_{M}^{*}$. The geometric series $\overline{S}_{P}(\lambda)$ can be calculated to be
$$
\overline{S}_{P}(\lambda)=(-1)^{|\Delta_{P}^{\lambda}|}e^{\langle \lambda,\, h_{P}  \rangle}\prod_{\alpha\in \Delta_{P}}\frac{1}{1-e^{-\langle\lambda,\,\alpha^{\vee}\rangle}},
$$
whence the result.

\end{proof}

Let $\re_{P}(\lambda)=\rd_{P}(\lambda)\overline{\rc}_{P}(\lambda)^{-1}$, they form a $(G,M)$-family, still denoted by $\re$. Lemma \ref{arthur lattice} expresses the lattice counting $\rw_{M}(h)$ as the volume of the product $(\rh\cdot \re)$ of $(G,M)$-families, where $\rh$ is the $(G,M)$-family associated to $h$. By the product formula (\ref{arthur prod}) of Arthur, we get
\begin{equation}\label{w mul 1}
\rw_{M}(h)=\sum_{Q\in \cf(M)}\rv_{M}^{Q}(h)\re_{Q}'.
\end{equation}

\begin{prop}\label{count in volume}
We have
$$
\rv_{\gamma}(x)=\sum_{Q=LN_{Q}\in \cf(M_{0})}\re_{Q}'\sum_{R=MN_{R}\in \cf^{L}(M_{0})} (-1)^{\dim(\ka_{M_{0}}^{M})}\cdot
\rv_{M_{0}}^{M}(\ec^{RN_{Q}}(x))\cdot \big[\ec^{Q}(x_{0})\big]'_{R}.
$$

\end{prop}

\begin{proof}
By lemma \ref{fund weight factor}, $\rv_{\gamma}(x)$ counts the number of lattice points in the polytope $\be(\lambda(x,x_{0}))$. Hence by the above formula,
$$
\rv_{\gamma}(x)=\rw_{M_{0}}(\lambda(x,x_{0}))=\sum_{Q\in \cf(M_{0})}\rv_{M_{0}}^{Q}(\lambda(x,x_{0}))\re_{Q}'.
$$
Now that $\lambda(x,x_{0})=\ec(x_{0})-\ec(x)$, for $Q=LN_{Q}\in \cf(M_{0})$, we have
\begin{align*}
\rv_{M_{0}}^{Q}(\lambda(x,x_{0}))&=\rv_{M_{0}}^{L}(\ec^{Q}(x_{0})-\ec^{Q}(x))
\\
&=\sum_{R\in \cf^{L}(M_{0})} 
\rv_{M_{0}}^{R}(-\ec^{Q}(x))\cdot \big[\ec^{Q}(x_{0})\big]'_{R}
\\
&=\sum_{R=MN_{R}\in \cf^{L}(M_{0})} (-1)^{\dim(\ka_{M_{0}}^{M})}\cdot
\rv_{M_{0}}^{M}(\ec^{RN_{Q}}(x))\cdot \big[\ec^{Q}(x_{0})\big]'_{R},
\end{align*}
where the last equality follows from a simple change of variable $\lambda\mapsto -\lambda$ in the definition (\ref{def gm volume}) of the volume of a $(G,M)$-family. Putting everything together, we get the assertion.

\end{proof}

\begin{thm}\label{f to w}
We have
$$
|F_{\gamma}(\fq)|=\vol(T^{1})\,|D(\gamma)|^{-\frac{1}{2}}\sum_{\substack{M,\,L\in \cl(M_{0})\\ M\subset L}}
(-1)^{\dim(\ka_{M_{0}}^{M})}J_{M_{0}}^{M}(\gamma, \mathbbm{1}_{\km\cap\kk})\,\rv_{M}^{L}(\ec(x_{0}))\, \re_{L}.
$$

\end{thm}

\begin{proof}

By lemma \ref{fund weight orb} and proposition \ref{count in volume}, we have
\begin{align*}
|F_{\gamma}(\fq)|&=\vol(T^{1})\int_{T\backslash G} \mathbbm{1}_{\kk}(\Ad(g^{-1})\gamma )\rv_{\gamma}(g)d\bar{g}
\\
&=\vol(T^{1})\sum_{Q=LN_{Q}\in \cf(M_{0})}\re_{Q}'
\sum_{R=MN_{R}\in \cf^{L}(M_{0})}
(-1)^{\dim(\ka_{M_{0}}^{M})}
\cdot \big[\ec^{Q}(x_{0})\big]'_{R}\cdot
\\
&\hspace{4cm}\int_{T\backslash G} \mathbbm{1}_{\kk}(\Ad(g^{-1})\gamma) \rv_{M_{0}}^{M}(\ec^{RN_{Q}}(g))\, d\bar{g}.
\end{align*}
With the Iwasawa decomposition $g=mnk$, $m\in M,n\in N_{R}N_{Q},k\in K$, the factor $\rv_{M_{0}}^{M}(\ec^{RN_{Q}}(g))$ simplifies to $\rv_{M_{0}}^{M}(m)$, and 
\begin{align*}
&\int_{T\backslash G} \mathbbm{1}_{\kk}(\Ad(g^{-1})\gamma)\rv_{M_{0}}^{M}(\ec^{RN_{Q}}(g))\,d\bar{g}
\\
&=\int_{T\backslash M}\int_{N_{R}N_{Q}}\int_{K} 
\mathbbm{1}_{\kk}(\Ad(k^{-1}n^{-1}m^{-1})\gamma)\rv_{M_{0}}^{M}(m)\,dkdnd\bar{m}
\\
&=|\det(\ad(\gamma);\,\kg/\km)|^{-\frac{1}{2}}\int_{T\backslash M} 
\mathbbm{1}_{\km\cap \kk}(\Ad(m^{-1})\gamma  )\rv_{M_{0}}^{M}(m)\,d\bar{m}
\\
&=|D(\gamma)|^{-\frac{1}{2}}J_{M_{0}}^{M}(\gamma,\mathbbm{1}_{\km\cap \kk}). 
\end{align*}
Plug in the above equation, $|F_{\gamma}(\fq)|$ equals
\begin{align*}
&\vol(T^{1})\,|D(\gamma)|^{-\frac{1}{2}}\sum_{Q=LN_{Q}\in \cf(M_{0})}\re_{Q}'
\sum_{R=MN_{R}\in \cf^{L}(M_{0})}
(-1)^{\dim(\ka_{M_{0}}^{M})}
\cdot \big[\ec^{Q}(x_{0})\big]'_{R}
\cdot
J_{M_{0}}^{M}(\gamma,\mathbbm{1}_{\km\cap \kk})
\\
&=\vol(T^{1})\,|D(\gamma)|^{-\frac{1}{2}}\sum_{L\in \cl(M_{0})}\re_{L}
\sum_{M\in \cl^{L}(M_{0})}
(-1)^{\dim(\ka_{M_{0}}^{M})}
\cdot 
 \ec_{M}^{L}(x_{0})\cdot
J_{M_{0}}^{M}(\gamma,\mathbbm{1}_{\km\cap \kk}),
\end{align*}
where the second equality follows the fact that $\ec^{Q}(x_{0})=\ec^{L}(x_{0})$ and that
$$
\sum_{R\in \cp(M)}\big[\ec^{L}(x_{0})\big]'_{R}=\ec_{M}^{L}(x_{0}),\quad 
\sum_{Q\in \cp(L)} \re_{Q}'=\re_{L}, 
$$
as implied by the equation (\ref{vol der}).

\end{proof}

The same technique applies to counting points on the {truncated} affine Springer fibers. 
Let $\Pi=(\mu_{P})_{P\in \cp(M_{0})}$ be a positive $(G,M_{0})$-orthogonal set with $\mu_{P}\in X_{*}(M_{0})$ for all $P\in \cp(M_{0})$. Let $\nu\in \Lambda_{G}$ be the image of $\mu_{P}$ under the projection $X_{*}(M_{0})=\Lambda_{M_{0}}\to \Lambda_{G}$, consider the truncated affine Springer fiber
$$
\xx_{\gamma}^{\nu}(\Pi)=\big\{x\in \xx_{\gamma}^{\nu}\mid \ec_{M_{0}}(x)\subset \Pi\big\}.
$$

\begin{thm}\label{f to w pi}
Assume that $\ec(x_{0})\subset \Pi$, then
{\small
$$
|\xx_{\gamma}^{\nu}(\Pi)(\fq)|=\vol(T^{1})\,|D(\gamma)|^{-\frac{1}{2}}
\sum_{Q=LN_{Q}\in \cf(M_{0})}\re_{Q}'\sum_{R=MN_{R}\in \cf^{L}(M_{0})} (-1)^{\dim(\ka_{M_{0}}^{M})}
J_{M_{0}}^{M}(\gamma, \mathbbm{1}_{\km\cap\kk}) (\Pi^{Q})'_{R}.
$$
}

\end{thm}

\begin{proof}

Observe that in the proof of lemma \ref{fund weight orb} and \ref{fund weight factor}, we use only the fact that $\ec(x)\subset \ec(x_{0})$ up to translation by $X_{*}(M_{0})$. Now that $\ec(x_{0})\subset \Pi$, we can express $|\xx_{\gamma}^{\nu}(\Pi)(\fq)|$ as a similar weighted orbital integral with the role of $\ec(x_{0})$ played by $\Pi$. We can then relate the weight factor to Arthur's weight factor in the same way to get the theorem.  

\end{proof}

We can reverse the above process to get an expression of $J_{M_{0}}(\gamma,\mathbbm{1}_{\kk})$ in terms of $|F_{\gamma}^{M}(\fq)|$, $M\in \cl(M_{0})$. Lemma \ref{arthur lattice} expresses $\rw_{M}(h)$ as the volume of the product $(\rh\cdot \re)$. Since $\rh=(\rh\cdot \re)\cdot \re^{-1}$, we get inversely
\begin{equation} \label{v in w}
\rv_{M_{0}}(h)=\sum_{Q\in \cf(M_{0})}
\rw_{M_{0}}^{Q}(h)(\re^{-1})'_{Q}.
\end{equation}

\begin{prop}
For any $x\in \xx_{\gamma}$, we have
$$
\rv_{M_{0}}(x)=\sum_{Q=LN_{Q}\in \cf(M_{0})}
(-1)^{\dim(\ka_{M_{0}}^{L})} \ec(x_{0})'_{Q}\sum_{R=MN_{R}\in \cf^{L}(M_{0})}
\rv_{\gamma}^{M}(f_{RN_{Q}}(x))\big[(\re^{-1})^{L}\big]'_{R}.
$$

\end{prop}

\begin{proof}

By construction, $\lambda(x,x_{0})=\ec(x_{0})-\ec(x)$, hence
\begin{align*}
\rv_{M_{0}}(\ec(x))&=\rv_{M_{0}}(\ec(x_{0})-\lambda(x,x_{0}))
=\sum_{Q=LN_{Q}\in \cf(M_{0})}
\rv_{M_{0}}^{Q}(-\lambda(x,x_{0}))\cdot \ec(x_{0})'_{Q}
\\
&=\sum_{Q=LN_{Q}\in \cf(M_{0})}
(-1)^{\dim(\ka_{M_{0}}^{L})}\rv_{M_{0}}^{L}(\lambda^{Q}(x,x_{0}))\cdot \ec(x_{0})'_{Q}
\\
&=\sum_{Q=LN_{Q}\in \cf(M_{0})}
(-1)^{\dim(\ka_{M_{0}}^{L})} \ec(x_{0})'_{Q}
\sum_{R=MN_{R}\in \cf^{L}(M_{0})}\rw_{M_{0}}^{R}(\lambda^{Q}(x,x_{0}))\big[(\re^{-1})^{L}\big]'_{R}.
\end{align*}
By definition, $\rw_{M_{0}}^{R}(\lambda^{Q}(x,x_{0}))$ counts the number of lattice points in the $R$-facet of $\lambda^{Q}(x,x_{0})$, hence it equals $\rv_{\gamma}^{M}(f_{RN_{Q}}(x))$, whence the assertion.

\end{proof}

\begin{thm}\label{w to f}
We have 
$$
J_{M_{0}}(\gamma,\mathbbm{1}_{\kk})=\vol(T^{1})^{-1}\sum_{\substack{M,L\in \cl(M_{0})\\M\subset L}}
(-1)^{\dim(\ka_{M_{0}}^{L})} |D^{M}(\gamma)|^{\frac{1}{2}}\cdot |F_{\gamma}^{M}(\fq)|\cdot (\re^{-1})^{L}_{M}\cdot\ec(x_{0})_{L}.
$$

\end{thm}

\begin{proof}

By the above result, $J_{M_{0}}(\gamma,\mathbbm{1}_{\kk})$ equals
\begin{align*}
&|D(\gamma)|^{\frac{1}{2}}\int_{T\backslash G}
\mathbbm{1}_{\kk}(\Ad(g^{-1})\gamma)\sum_{Q=LN_{Q}\in \cf(M_{0})}
(-1)^{\dim(\ka_{M_{0}}^{L})} \ec(x_{0})'_{Q}\cdot
\\
&\hspace{4cm}\sum_{R=MN_{R}\in \cf^{L}(M_{0})}
\rv_{\gamma}^{M}(f_{RN_{Q}}(g))\big[(\re^{-1})^{L}\big]'_{R}d\bar{g}.
\end{align*}
With the Iwasawa decomposition, we calculate
\begin{align*}
&|D(\gamma)|^{\frac{1}{2}}\int_{T\backslash G}
\mathbbm{1}_{\kk}(\Ad(g^{-1})\gamma)\rv_{\gamma}^{M}(f_{RN_{Q}}(g))d\bar{g}
\\
&=|D(\gamma)|^{\frac{1}{2}}\int_{T\backslash M}\int_{N_{R}N_{Q}}\int_{K}
\mathbbm{1}_{\kk}(\Ad(k^{-1}n^{-1}m^{-1}) \gamma)\rv_{\gamma}^{M}(m)\,dkdnd\bar{m}
\\
&=|D^{M}(\gamma)|^{\frac{1}{2}}
\int_{T\backslash M}
\mathbbm{1}_{\km\cap \kk}(\Ad(m^{-1}) \gamma )\rv_{\gamma}^{M}(m)\,d\bar{m}
\\
&=|D^{M}(\gamma)|^{\frac{1}{2}}\cdot \vol(T^{1})^{-1}\cdot|F_{\gamma}^{M}(\fq)|,
\end{align*}
Plug in the above equation, $J_{M_{0}}(\gamma,\mathbbm{1}_{K})$ equals 
\begin{align*}
&\vol(T^{1})^{-1}\sum_{Q=LN_{Q}\in \cf(M_{0})}
(-1)^{\dim(\ka_{M_{0}}^{L})} \ec(x_{0})'_{Q}\cdot
\sum_{R=MN_{R}\in \cf^{L}(M_{0})}
|D^{M}(\gamma)|^{\frac{1}{2}}\cdot |F_{\gamma}^{M}(\fq)|\cdot\big[(\re^{-1})^{L}\big]'_{R}
\\
&=\vol(T^{1})^{-1}\sum_{L\in \cl(M_{0})}
(-1)^{\dim(\ka_{M_{0}}^{L})} \ec(x_{0})_{L}\cdot
\sum_{M\in \cl^{L}(M_{0})}
|D^{M}(\gamma)|^{\frac{1}{2}}\cdot |F_{\gamma}^{M}(\fq)|\cdot (\re^{-1})^{L}_{M},\end{align*}
where in the second equation we have used the equation (\ref{vol der}) to get
$$
\sum_{Q\in \cp(L)}\ec(x_{0})'_{Q}=\ec(x_{0})_{L},\quad \sum_{R\in \cp^{L}(M)}\big[(\re^{-1})^{L}\big]'_{R}=(\re^{-1})^{L}_{M}.
$$
Rearrange the terms, we get the assertion.

\end{proof}

In \cite{chen2}, the transition between weighted orbital integrals and counting points on the fundamental domains was extended to general Levi subgroup $M\in \cl(M_{0})$. This can be achieved also with our current method. Indeed, we have a reduction formula for $J_{M}(\gamma,\varphi)$.

\begin{prop}\label{weight orb reduction}

Let $\gamma\in \kg$ be an element satisfying $M_{\gamma}=G_{\gamma}$, let $\varphi\in \cc_{c}^{\infty}(\kg)$, then 
$$
J_{M}(\gamma,\varphi)=\sum_{L\in \cl(M_{0})}
\theta_{M_{0}}^{G}(M,L)\cdot J_{M_{0}}^{L}(\gamma,\varphi^{Q_{L}^{\xi}}),\quad M\in \cl(M_{0}),
$$
where $Q_{L}^{\xi}\in \cp(L)$ is uniquely determined by a generic element $\xi\in \ka_{M}^{L}$ as in equation $(\ref{proj vol})$.

\end{prop}

\begin{proof}

By equation (\ref{proj vol}), we have
\begin{align*}
J_{M}(\gamma,\varphi)
&=|D(\gamma)|^{\frac{1}{2}} \int_{T\backslash G}
\varphi(\Ad(g^{-1})\gamma) 
\sum_{L\in \cl(M_{0})}\theta_{M_{0}}^{G}(M,L) \rv_{M_{0}}^{Q_{L}^{\xi}}(g)
\,d\bar{g}
\\
&=|D(\gamma)|^{\frac{1}{2}}
\sum_{L\in \cl(M_{0})}\theta_{M_{0}}^{G}(M,L)
 \int_{T\backslash L} \int_{N_{Q_{L}^{\xi}}} \int_{K}
\varphi(\Ad(k^{-1}n^{-1}l^{-1})\gamma ) 
 \rv_{M_{0}}^{L}(l)
\,d\bar{l}
\\
&=
\sum_{L\in \cl(M_{0})}\theta_{M_{0}}^{G}(M,L)
 |D^{L}(\gamma)|^{\frac{1}{2}}
 \int_{T\backslash L} 
\varphi^{Q_{L}^{\xi}}(\Ad(l^{-1})\gamma  ) 
 \rv_{M_{0}}^{L}(l)
\,d\bar{l}.
\end{align*}

\end{proof}

\begin{rem}

The reduction formula is reminiscent of the corollary 3.10 in \cite{chen2}, which reduces the counting points of the intermediate fundamental domain to that of the fundamental domains.  

\end{rem}

\section{Rationality from homogeneity}\label{r from h}

This section is devoted to the proof of theorem \ref{main}. We work with the group $\bg=\gl_{d+1}, d\ge 1$, let $\ba$ be its maximal torus of the diagonal matrices, let $\gamma\in \ka$ be a regular integral element. Following Goresky, Kottwitz and MacPherson \cite{gkm4}, we define the \emph{root valuation datum} of $\gamma$ to be the function
$$
R_{\gamma}:\Phi=\Phi(\bg,\ba)\to \bz,\quad \alpha\mapsto \val(\alpha(\gamma)).
$$

\begin{lem}

There exists a set $\Delta=\{\alpha'_{i}\}_{i=1}^{d}$ of simple roots of $\Phi$ such that for any root $\alpha=\sum_{i=1}^{d}a_{i}\alpha_{i}'\in \Phi$, we have
$$
R_{\gamma}(\alpha)=\min\vparen{R_{\gamma}(\alpha'_{i})\mid a_{i}\neq 0}.
$$ 
In this case, we say that $\gamma$ is \emph{in minimal form} with respect to $\Delta$. 
\end{lem}

\begin{proof}

The question involves only the root system, we proceed by induction on the dimension of the root system. The assertion holds obviously in dimension $1$. Suppose that it holds for  dimension less than $d$.

For $m\in \bz$, let $\Phi_{\gamma,m}=\vparen{\alpha\in \Phi \mid R_{\gamma}(\alpha)\ge m}$, then we get a decreasing filtration
\begin{equation}\label{filt roots}
\Phi=\Phi_{\gamma,0}\supset \Phi_{\gamma,1}\supset \Phi_{\gamma,2}\supset \cdots.
\end{equation}
According to Goresky, Kottwitz and MacPherson \cite{gkm4}, lemma 4.8.2, $\Phi_{\gamma,m}$ is $\bq$-closed, i.e. any $\bq$-linear combination of roots in $\Phi_{\gamma,m}$, if it is a root, must be an element of $\Phi_{\gamma,m}$. In particular, $\Phi_{\gamma,m}$ is itself a root system. Let $m_{0}$ be the smallest integer such that $\Phi_{\gamma,0}\supsetneq \Phi_{\gamma,m}$. By induction hypothesis, there exists a set $\Delta'$ of simple roots in $\Phi_{\gamma,m_{0}}$ such that the restriction of the root valuation datum to $\Phi_{\gamma,m_{0}}$ satisfies the conclusion of the assertion. By \cite{bourbaki lie}, chap. 6, \S7, prop. 24, the set $\Delta'$ can be extended to a set of simple roots $\Delta$ of $\Phi$. It is easy to verify that $\Delta$ satisfies the property in the assertion with the strong triangular inequality of the valuation.

\end{proof}

Since the affine Springer fiber $\xx_{\gamma}$ depends only on the conjugacy class of $\gamma$ and that all the sets of simple roots are conjugate under the Weyl group $W$, by conjugating $\gamma$ with a suitable element in $W$, we can assume that $\gamma$ is in minimal form with respect to the set of simple roots $\Delta_{0}=\{\alpha_{i}\}_{i=1}^{d}$ associated to the standard Borel subgroup $\mathbf{B}_{0}$ of the upper triangular matrices, and we use also the name \emph{root valuation datum} of $\gamma$ for the $d$-tuple $\bnn=\big(R_{\gamma}(\alpha_{i})\big)_{i=1}^{d}$. Obviously, for any $\bnn\in \bz_{\ge 0}^{d}$, we can find $\gamma$ in minimal form with root valuation datum $\bnn$. 
Moreover, as $\xx_{\gamma}=\xx_{\gamma+a}$ for any $a\in \ka$, we can always assume that the minimum of $R_{\gamma}(\alpha_{i}), i=1,\cdots, d$, equals the minimum of the valuations of the eigenvalues of $\gamma$.
These notations and conventions generalize to the Levi subgroups $M\in \cl(A)$ in the obvious way, we don't write it down.

\subsection{Generating series of weighted orbital integrals}

To begin with, we investigate a similar generating series of the weighted orbital integrals. Suppose that $\gamma$ is in minimal form with respect to $\Delta_{0}$ and let $\bnn$ be its root valuation. Let $m_{1}<\cdots<m_{l}$ be the breaking points of the filtration (\ref{filt roots}), i.e. the filtration reduces to
\begin{equation}\label{filt roots n}
\Phi_{\gamma,m_{1}}\supsetneq \Phi_{\gamma,m_{2}}\supsetneq \cdots \supsetneq \Phi_{\gamma,m_{l}}.
\end{equation}
Then for any root $\alpha\in \Phi_{\gamma, m_{i}}\backslash \Phi_{\gamma, m_{i+1}}$, we have $R_{\gamma}(\alpha)=m_{i}$. Let ${M}_{(i)}\in \cl(A)$ be the unique Levi subgroup with root system $\Phi_{\gamma, m_{i}}$. 

We have the following basic reduction: Let $\gamma'={\ep^{-m_{1}}}\gamma$, then $\gamma'$ remains integral with root valuation datum dropped by $m_{1}$. In particular, for the roots in  $\kg/\km_{(2)}$, the root valuations are all $0$, and we can write
$$
\gamma'=s+\eta, 
$$ 
with $s\in \ka_{M_{(2)}}^{\reg}$ and $\eta\in \km_{(2)}\cap U_{s}^{+}$ both integral and that the eigenvalues of $\ad(s)$ on $\kg/\km_{(2)}$ all have valuation $0$. With the Arthur-Shalika expansion, i.e. theorem \ref{shalika}, for $f\in \cc_{c}(\kg/\kk)$, we can reduce the study of the weighted orbital integral $J_{A}(\gamma,f)$ to that of the Arthur-Shalika germs $g_{A}^{L}(\gamma,u)$ and the unipotent weighted orbital integrals $J_{L}(u,f)$, for $L\in \cl(A)$ and $u\in \nn_{L}(F)$. 
By theorem \ref{shalika homogeneity}, $g_{A}^{L}(\gamma,u)$ can be expressed in terms of $g_{A}^{M}(\gamma',v)=g_{A}^{M}(s+\eta,v)$ with $M\in \cl^{L}(A)$ and $v\in \nn_{M}(F)$.
We can make a further parabolic reduction for $g_{A}^{M}(s+\eta,v)$: with theorem \ref{shalika} and theorem \ref{germ descent}, we get two expansions of $J_{A}(s+\eta,f)$, which give rise to the linear equation 
\begin{equation}\label{germ hyper reduce}
\sum_{M\in \cl(A)}\sum_{v\in [\nn_{M}(F)]}
g_{A}^{M}(s+\eta, v) J_{M}(v,f)
=
\sum_{M'\in \cl^{M_{(2)}}(A)}\sum_{u'\in [\nn_{M'}(F)]}
g_{A}^{M'}(\eta,u')J_{M'}(s+u',f).  
\end{equation}
Recall that $g_{A}^{M}$ is constructed by induction on $M$, hence it is enough to work out $g_{A}^{G}$. 
Let $f$ run through a dual basis of the distributions $J_{G}(u,\,\cdot\,)$, $u\in [\nn_{G}(F)]$, we get a linear system.
Solving them inductively in $M$ and $v\in [\nn_{M}(F)]$, we can express $g_{A}^{M}(s+\eta,v)$ in terms of $g_{A}^{M'}(\eta,u')$ with coefficients in $J_{M}(v,f)$ and $J_{M'}(s+u',f)$. 
Hence the problem is further reduced to the study of $g_{A}^{M'}, M'\in \cl^{M_{(2)}}(A)$, i.e. to the study of the local behvior of $J_{A}^{M_{(2)}}(\eta,\,\cdot\,)$ for the group $M_{(2)}$. Briefly, the study of $J_{A}(\gamma,\,\cdot\,)$ can be reduced to that of   
$$
\begin{cases}
J_{M}(u,\,\cdot\,),& \text{ for }M\in \cl(A), u\in [\nn_{M}(F)],\\
J_{M'}(s+u',\,\cdot\,),& \text{ for }M'\in \cl^{M_{(2)}}(A),  u'\in [\nn_{M'}(F)],\\
J_{A}^{M_{(2)}}(\eta,\,\cdot\,).
\end{cases}
$$
This process can be carried on for $J_{A}^{M_{(2)}}(\eta,\,\cdot\,)$ and so on. At the end, the study of $J_{A}(\gamma,\,\cdot\,)$ can be reduced to 
\begin{equation}\label{finite uni}
\begin{cases}
J_{M}^{M_{(i)}}(u,\,\cdot\,),& \text{ for }M\in \cl^{M_{(i)}}(A), u\in [\nn_{M}(F)],\\
J_{M'}^{M_{(i)}}(s_{i}+u',\,\cdot\,),& \text{ for }
M'\in \cl^{M_{(i+1)}}(A),  u'\in [\nn_{M'}(F)],\\
\end{cases}
\end{equation}
where $s_{i}\in \ka_{M_{(i+1)}}^{\reg}$ are integral such that the eigenvalues of $\ad(s_{i})$ on $\km_{(i)}/\km_{(i+1)}$ all have valuation $0$, $i=1,\cdots,l$. For $L\in \cl(A)$, let $\ka_{L}^{\reg,0}$ be the open subset of $\ka_{L}^{\reg}$ consisting of the integral elements $s$ such that the eigenvalues of $\ad(s)$ on $\kg/\kl$ all have valuation $0$.

\begin{lem}\label{key lemma}

Let $L\in \cl(A)$ be such that $\dim(A_{L})=\dim(A_{G})+1$, let $s\in \ka_{L}^{\reg,0}$. Then for any $f\in \cc_{c}(\kg/\kk)$, $M\in \cl^{L}(A)$ and $u\in [\nn_{M}(F)]$, as a function of $s\in \ka_{L}^{\reg,0}$, the weighted  orbital integral $J_{M}(s+u,f)$ is constant.

\end{lem}

\begin{proof}

Our original proof has been computational and clumsy. Prof. Jean-Loup Waldspurger very kindly offers us a short and elegant proof, we reproduce it here. 

Let $s=s_{G}+s^{G}$ under the decomposition $\ka_{A}=\ka_{G}\oplus \ka_{A}^{G}$. To begin with, the integral $J_{M}(s+u,f)$ doesn't depend on the component $s_{G}\in \ka_{G}$. Indeed, without loss of generality, we can suppose that $f$ is a tensor product $f_{G}\otimes f^{G}$ with $f_{G}\in \cc_{c}(\ka_{G}/(\ka_{G}\cap\kk))$ and $f^{G}\in \cc_{c}(\kg^{\der}/(\kg^{\der}\cap \kk))$. Then
$$
J_{M}(s+u,f)=J^{A_{G}}(s_{G},f_{G})\cdot J_{M\cap G^{\der}}^{G^{\der}}(s^{G}+u,f^{G}),
$$
and $J^{A_{G}}(s_{G},f_{G})=f_{G}(0)$ is independent of $s_{G}$. 

We can now restrict to $s\in \ka_{L}^{G}\cap \ka_{L}^{\reg,0}$. As $\ka_{L}^{G}$ is one dimensional, the assertion can be reformulated as the constance of the function $t\mapsto J_{M}(ts+u,f)$ for $t\in \co_{F}^{\times}$. Recall that the integral $J_{M}(\gamma,f)$ depends only on the $M(F)$-conjugacy class of $\gamma$, we get
$$
J_{M}(ts+u,f)=J_{M}(ts+tu,f)=J_{M}(s+u,f^{t})
$$
as $ts+u$ and $ts+tu$ are $M(F)$-conjugate. The question is then reduced to showing the constance of the function $t\mapsto J_{M}(s+u,f^{t})$ for $t\in \co_{F}^{\times}$. More generally, we can consider the function $t\mapsto J_{M}(\gamma,f^{t})$ for an integral element $\gamma\in \km$ and $t\in \co_{F}^{\times}$. By theorem \ref{shalika}, we have
$$
J_{M}(\gamma,f^{t})=\sum_{L\in \cl(M)}\sum_{\nu\in [\nn_{L}(F)]} g_{M}^{L}(\gamma,\nu) J_{L}(\nu, f^{t}).
$$
But $J_{L}(\nu, f^{t})=J_{L}(t\nu,f)=J_{L}(\nu,f)$ is clearly independent of $t$, as $\nu$ and $t\nu$ are $L(F)$-conjugate. This proves the constance of the function that we start with.

\end{proof}

\begin{thm}\label{main weighted}

Let $\bg=\gl_{d+1}$, let $\gamma\in \ka$ be a regular integral element, then for any $f\in \cc_{c}(\kg/\kk)$, the weighted orbital integral $J_{A}(\gamma,\,f)$ depends only on the root valuation datum $\bnn$ of $\gamma$, and we denote it by $J_{A}(\bnn, f)$. Moreover, the generating series
$$
J_{f}(\mathbf{t}):=\sum_{n_{1}=0}^{+\infty}\cdots \sum_{n_{d}=0}^{+\infty}  J_{A}\big((n_{1},\cdots,n_{d}), f\big) t_{1}^{n_{1}}\cdots t_{d}^{n_{d}}\in \bc [\![ t_{1},\cdots,t_{d}]\!]
$$
is a rational fraction, \textup{i.e.} it is an element of $\bc( t_{1},\cdots,t_{d})$, and the denominator is divided by
$$
\prod_{i=1}^{d}(1-t_{i}) \cdot
\prod_{\substack{M\in \cl(A)\\ M=M_{1}\times\cdots \times M_{r}}} 
\prod_{M'\in \cl^{M}(A)}\prod_{\substack{v\in [\cn_{M'}(F)]\\ v=(v_{1},\cdots, v_{r})}} \prod_{j=1}^{r}\Big(1-q^{-d^{M'\cap M_{j}}(v_{j})}\prod_{\alpha_{i}\in \Delta_{M_{j}}}t_{i}\Big), 
$$
where $M=M_{1}\times \cdots\times M_{r}$ is the factorizations into linear groups, and $v=(v_{1}, \cdots, v_{r})$ is the coordinate with respect to the corresponding decomposition of $\km$.

\end{thm}

By theorem \ref{shalika}, the above theorem is an immediate consequence of the following theorem \ref{shalika generating series} about the Arthur-Shalika germs, we state it separately because it has independent interest.


\begin{thm}\label{shalika generating series}

Let $\bg=\gl_{d+1}$, let $\gamma\in \ka$ be a regular integral element, let $u\in \cn_{G}(F)$. For the Levi subalgebra $\km$ containing $u$, the Arthur-Shalika germ $g_{A}^{M}(\gamma,\,u)$ depends only on the root valuation datum $\bnn$ of $\gamma$, and we denote it by  $g_{A}^{M}(\bnn, u)$. The generating series
$$
g_{A}^{M}(u;\mathbf{t}):=\sum_{\bnn\in \bz_{\ge 0}^{d}}  g_{A}^{M}\big(\bnn, u\big) t_{1}^{n_{1}}\cdots t_{d}^{n_{d}}
$$
is a rational fraction, and the denominator is divided by
\begin{align*}
&\prod_{i=1}^{d}(1-t_{i}) \cdot
\prod_{j=1}^{r} \Bigg[\Big(1-q^{-d^{M_{j}}(u_{j})}\prod_{\alpha_{i}\in \Delta_{M_{j}}}t_{i}\Big) 
\prod_{\substack{J\in \cl^{M}(A)\\ J\neq M}}  \prod_{v\in [\cn_{J}(F)]} \Big(1-q^{-d^{J\cap M_{j}}(v_{j})}\prod_{\alpha_{i}\in \Delta_{M_{j}}}t_{i}\Big)\Bigg]\cdot
\\
&\prod_{\substack{M'\in \cl^{M}(A),\, M'\neq M\\ M'=M_{1}'\times\cdots \times M_{r'}'}} 
\prod_{M''\in \cl^{M'}(A)}\prod_{\substack{v'\in [\cn_{M''}(F)]\\ v'=(v'_{1},\cdots, v'_{r'})}} \prod_{j=1}^{r'}\Big(1-q^{-d^{M''\cap M'_{j}}(v'_{j})}\prod_{\alpha_{i}\in \Delta_{M'_{j}}}t_{i}\Big), 
\end{align*}
where $M=M_{1}\times \cdots\times M_{r}$ and $M'=M_{1}'\times\cdots \times M_{r'}'$ are factorizations into linear groups, and $u=(u_{1}, \cdots, u_{r})$, $v=(v_{1}, \cdots, v_{r})$ and $v'=(v'_{1},\cdots, v'_{r'})$ are the coordinates with respect to the corresponding decomposition of $\km$ and $\km'$ respectively.

\end{thm}

\begin{proof}

To begin with, note that the Arthur-Shalika germ $g_{A}^{M}(\gamma, u)$ has a factorization property:
Suppose that $M$ admits a factorization into linear groups $M=M_{1}\times\cdots\times M_{r}$, then we can write $A=A_{1}\times\cdots\times A_{r}$ accordingly, and with respect to the associated decomposition of $\km$ we can write $\gamma=(\gamma_{1},\cdots, \gamma_{r})$ and $u=(u_{1},\cdots, u_{r})$, then
\begin{equation}\label{factorize g}
g_{A}^{M}(\gamma, u)=g_{A_{1}}^{M_{1}}(\gamma_{1}, u_{1})\times \cdots \times g_{A_{r}}^{M_{r}}(\gamma_{r}, u_{r}).
\end{equation}
Indeed, for a factorizable test function $\varphi=\varphi_{1}\otimes\cdots\otimes \varphi_{r}\in \cc_{c}^{\infty}(\km_{1})\times \cdots\times \cc_{c}^{\infty}(\km_{r})$, we have 
$$
J_{A}^{M}(\gamma, \varphi)=J_{A_{1}}^{M_{1}}(\gamma_{1}, \varphi_{1})\times\cdots\times  J_{A_{r}}^{M_{r}}(\gamma_{r}, \varphi_{r}),
$$
and the factorization property (\ref{factorize g}) follows from it. 


For the dependence of $g_{A}^{M}(\gamma, u)$ on the root valuation datum of $\gamma$, it is enough to prove it for $M=G$ as we have explained above. 
We apply the basic reduction at the beginning of the section and we use the notations over there. 
With the homogeneity of $g_{A}^{G}(\gamma, u)$, we can reduce to the case that  $\gamma'=s+\eta$ with $s\in \ka_{M_{(2)}}^{\reg,0}$ and $\eta\in \km_{(2)}\cap U_{s}^{+}$ both integral. By theorem \ref{shalika homogeneity}, this process depends only on the root valuation datum of $\gamma$. 
By the transitivity of the parabolic descent, we can assume that $\dim(A_{M_{(2)}})=\dim(A_{G})+1$. As before, we get two expansions of $J_{A}(s+\eta,f)$ as equation (\ref{germ hyper reduce}). 
Let $f$ run through a dual basis of the linear functionals $J_{G}(u, \cdot)$, $u\in [\nn_{G}(F)]$ in the equation (\ref{germ hyper reduce}), we get then a non-degenerate linear system for the germs $g_{A}^{G}(s+\eta, u)$.  Solving it, we get a a linear expression of $g_{A}^{G}(s+\eta, u)$ in terms of 
$$g_{A}^{M}(s+\eta, v), J_{M}(v, f), \quad M\in \cl^{G}(A), M\neq G, v\in [\nn_{M}(F)],
$$ 
and 
$$
g_{A}^{M'}(\eta, v'), J_{M'}(s+v',f), \quad M'\in \cl^{M_{(2)}}(A), v'\in [\nn_{M'}(F)].
$$ 
Iterate the above process for $g_{A}^{M}(s+\eta, v)$, we get a linear expression of $g_{A}^{G}(s+\eta, u)$ in terms of 
$$
J_{M}(v, f), \quad M\in \cl^{G}(A), M\neq G, v\in [\nn_{M}(F)],
$$ 
and
$$
g_{A}^{M'}(\eta, v'), J_{M'}(s+v',f), \quad M'\in \cl^{M_{(2)}}(A), v'\in [\nn_{M'}(F)].
$$ 
By lemma \ref{key lemma}, $J_{M'}(s+v',f)$ is constant as a function of $s\in \ka_{M_{(2)}}^{\reg,0}$. We summarize the above discussions as

\begin{lem}\label{reduce gm levi 1}

Under the above setting, the Arthur-Shalika germ $g_{A}^{G}(s+\eta, u)$ can be expressed as a linear combination of 
$g_{A}^{M'}(\eta, v')$ with $M'\in \cl^{M_{(2)}}(A), v'\in [\nn_{M'}(F)]$.

\end{lem}

\noindent The dependence of $g_{A}^{G}(\gamma, u)$ on the root valuation datum of $\gamma$ then follows by induction on $\rk(\bfg)$.

Before proceeding, we make some comments on the power series $g_{A}^{M}(u;\btt)$. With the factorization property (\ref{factorize g}), it is clear that $g_{A}^{M}(\gamma, u)$ depends only on the root valuation datum of $\gamma$ with respect to $M$, i.e. the restriction of the root valuation datum $\bnn$ to the set $\Delta_{M}$. For this dependence property, we denote $\bar{g}_{A}^{M}(\bmm, u)=g_{A}^{M}(\gamma, u)$, and we should have considered the power series
$$
\bar{g}_{A}^{M}(u; \btt)=\sum_{\bmm\in \bz_{\ge 0}^{\Delta_{M}}} \bar{g}_{A}^{M}(\bmm, u)\prod_{\alpha_{i}\in \Delta_{M}}t_{i}^{m_{\alpha_{i}}}.
$$
Let $M=M_{1}\times \cdots M_{r}$ be the factorization into linear groups, let $u=(u_{1}, \cdots, u_{r})$ be the coordinate with respect to the decomposition $\km=\km_{1}\oplus \cdots \oplus \km_{r}$. Let $A_{i}=A\cap M_{i}$  and let $\bmm_{i}$ be the restriction of $\bmm$ to the set of simple roots $\Delta_{M_{i}}$ in $\Phi(M_{i}, A_{i})$, let $\btt_{i}$ be the set of indeterminants indexed by elements in $\Delta_{M_{i}}$, then we have the factorization
\begin{equation}\label{ga reduce 1}
\bar{g}_{A}^{M}(u; \btt)=\prod_{i=1}^{r}g_{A_{i}}^{M_{i}}(u_{i}, \btt_{i}).
\end{equation}
Moreover, we have the simple transformation
\begin{equation}\label{ga reduce 2}
g_{A}^{M}(u, \btt)=\bar{g}_{A}^{M}(u,\btt) \prod_{\alpha_{i}\in \Delta_{G}-\Delta_{M}}(1-t_{i})^{-1}. 
\end{equation}
With the property (\ref{ga reduce 1}) and (\ref{ga reduce 2}), it is enough to prove the remaining assertions for $M=G$.

Let $\Phi_{\bullet}(\bnn)$ be the reduced strict filtration associated to (\ref{filt roots n}) with $\Phi_{i}(\bnn)=\Phi_{\gamma,m_{i}}, i=1,\cdots, l$. We regroup the summation in the definition of $g_{A}^{M}(u;\mathbf{t})$ according to the filtration of roots $\Phi_{\bullet}(\bnn)$ as follows: Let $\Phi_{\bullet}$ be a filtration  
$$
\Phi_{1}=\Phi\supsetneq \Phi_{2}\supsetneq \cdots \supsetneq \Phi_{l}
$$ 
of $\Phi$ by $\bq$-closed sub root systems. We denote loosely $\Phi_{\bullet}(\bnn)=\Phi_{\bullet}$ if the filtrations have the same length and that $\Phi_{m_{i}}(\bnn)=\Phi_{i}$ for all the breaking points $m_{i}$ of the filtration $\Phi_{\bullet}(\bnn)$. 
For $i=1, \cdots, l$, let $M_{(i)}$ be the Levi subgroup with root system $\Phi_{i}$. 
Let
$$
S^{M}_{\Phi_{\bullet}}(u;\mathbf{t}):=\sum_{{\bnn\in \bz_{\ge 0}^{d},\, \Phi_{\bullet}(\bnn)=\Phi_{\bullet}}}  
g_{A}^{M}(\bnn, u) t_{1}^{n_{1}}\cdots t_{d}^{n_{d}}.
$$
Then the assertion is reduced to showing that $S^{M}_{\Phi_{\bullet}}(u;\mathbf{t})$ is rational, as there are only finitely many filtrations $\Phi_{\bullet}$.

We break the summation in $S^{G}_{\Phi_{\bullet}}(u;\mathbf{t})$ into two terms, $S_{1}(u;\mathbf{t})$ and $S_{2}(u;\mathbf{t})$, according to whether $m_{1}=0$ or not. 
For the term $S_{1}(u;\mathbf{t})$, suppose that $M_{(2)}$ has a factorization into linear groups $M_{(2)}=M_{2, 1}\times\cdots M_{2,r_{2}}$. Let $A_{i}=A\cap M_{2, i}$ and $\bnn_{i}$ be the restriction of the root valuation datum $\bnn$ to the simple roots in $\Delta_{M_{2,i}}$. 
By lemma \ref{reduce gm levi 1} and the factorization property (\ref{factorize g}), we have

\begin{lem}\label{reduce gm levi 2}

The power series $S_{1}(u;\btt)$ can be written as a linear combination of 
$$
\prod_{i=1}^{r_{2}} S_{\Phi_{\bullet}^{M_{2,i}}}^{M'\cap M_{2, i}}(v_{i}', \btt_{i}), \quad M'\in \cl^{ M_{(2)}}(A), v'\in [\nn_{M'}(F)],
$$
where $\Phi_{\bullet}^{M_{2,i}}$ is the reduced strict filtration on $\Phi(M_{2, i}, A_{i})$ induced by $\Phi_{\bullet}$, and $\btt_{i}$ is indeterminants indexed by $\Delta_{M_{2, i}}$.

\end{lem}

\noindent The rationality of $S_{1}(u;\mathbf{t})$ then follows by induction.

Now that the rationality of $S_{1}(u;\mathbf{t})$ is established, with the homogeneity of the Arthur-Shalika germs, we can make the equation
$$
S^{G}_{\Phi_{\bullet}}(u;\mathbf{t})=S_{1}(u;\mathbf{t})+S_{2}(u;\mathbf{t})
$$
into a recurrence equation. Let $\mathbf{1}=(1,\cdots,1)\in \bz^{d}$. By theorem \ref{shalika homogeneity}, we have
\begin{align*}
S^{G}_{\Phi_{\bullet}}(u;\mathbf{t})
&=
S_{1}(u;\mathbf{t})+\sum_{{\bnn\in \bn^{d},\, \Phi_{\bullet}(\bnn)=\Phi_{\bullet}}}  
g_{A}^{G}(\bnn, u) t_{1}^{n_{1}}\cdots t_{d}^{n_{d}}
\\
&=S_{1}(u;\mathbf{t})+\sum_{{\bnn\in \bz_{\ge 0}^{d},\, \Phi_{\bullet}(\bnn)=\Phi_{\bullet}}}  
g_{A}^{G}(\bnn+\mathbf{1}, u) t_{1}^{n_{1}+1}\cdots t_{d}^{n_{d}+1}
\\
&=S_{1}(u;\mathbf{t})+\sum_{{\bnn\in \bz_{\ge 0}^{d},\, \Phi_{\bullet}(\bnn)=\Phi_{\bullet}}}  
q^{-d^{G}(u)}\sum_{M\in \cl(A)}
\sum_{v\in [\nn_{M}(F)]} g_{A}^{M}(\bnn,v){\rc}_{M}^{G}(v,\ep)
\\
&\hspace{8cm} \cdot \big[v_{M}^{G}:u/\ep\big]
t_{1}^{n_{1}+1}\cdots t_{d}^{n_{d}+1}
\\
&=S_{1}(u;\mathbf{t})+  
q^{-d^{G}(u)} \Bigg[S^{G}_{\Phi_{\bullet}}(u/\ep;\mathbf{t})t_{1}\cdots t_{d}+
\\
&\hspace{3cm} 
\sum_{\substack{M\in \cl(A)\\ M\neq G}}
\sum_{v\in [\nn_{M}(F)]} 
\big[v_{M}^{G}:u/\ep\big] {\rc}_{M}^{G}(v,\ep)
S^{M}_{\Phi_{\bullet}}(v;\mathbf{t}) t_{1}\cdots t_{d}\Bigg],
\end{align*}
where for the last equality we have used the fact that ${\rc}_{G}^{G}(u/\ep,\ep) =1$.
Since $u$ and $u/\ep$ belongs to the same nilpotent orbit in $\nn_{G}(F)$, we have $\big[v_{M}^{G}:u/\ep\big]=\big[v_{M}^{G}:u\big]$ and
$$
S^{G}_{\Phi_{\bullet}}(u;\mathbf{t})=S^{G}_{\Phi_{\bullet}}(u/\ep;\mathbf{t}),
$$
and the above equation reduces to
\begin{align}
&\Big(1-q^{-d^{G}(u)}t_{1}\cdots t_{d} \Big)S^{G}_{\Phi_{\bullet}}(u;\mathbf{t}) \nonumber
\\
&=S_{1}(u;\mathbf{t})+  
q^{-d^{G}(u)}
\sum_{\substack{M\in \cl(A)\\ M\neq G}}
\sum_{v\in [\nn_{M}(F)]} 
\big[v_{M}^{G}:u\big] {\rc}_{M}^{G}(v,\ep)
S^{M}_{\Phi_{\bullet}}(v;\mathbf{t}) t_{1}\cdots t_{d}. \label{recursion}
\end{align}
The rationality of $S^{G}_{\Phi_{\bullet}}(u;\mathbf{t})$ then follows by induction.

By the recursion (\ref{recursion}), we obtain that the denominator of $S_{\Phi_{\bullet}}^{G}(u,\btt)$ divides $\big(1-q^{-d^{G}(u)}t_{1}\cdots t_{d} \big)$ times that of $S_{1}(u;\mathbf{t})$ and $S^{M}_{\Phi_{\bullet}}(v;\mathbf{t})$. 
By lemma \ref{reduce gm levi 2}, the denominator of $S_{1}(u;\mathbf{t})$ can be further reduced. 
At the end, we arrive at $S^{A}_{\Phi_{\bullet}}(0, \btt)$, which has denominator $\prod_{i=1}^{d}(1-t_{i})$. 
For $a=1,\cdots, l$, let $M_{(a)}=M_{a, 1}\times\cdots \times M_{a, r_{a}}$ be the factorization into linear groups, we conclude that the denominator of $S_{\Phi_{\bullet}}^{G}(u,\btt)$ divides
\begin{align*}
&\prod_{i=1}^{d}(1-t_{i})\cdot \big(1-q^{-d^{G}(u)}t_{1}\cdots t_{d} \big)\prod_{\substack{M\in \cl(A)\\ M\neq G}}\prod_{v\in [\nn_{M}(F)]} \big(1-q^{-d^{M}(v)}t_{1}\cdots t_{d} \big)
\\
\\
&\prod_{a=1}^{l} \prod_{M'\in \cl^{M_{(a)}}(A)}\prod_{\substack{u'\in [\nn_{M'}(F)]\\ u'=(u'_{1},\cdots, u'_{r_{a}})}}
\prod_{j=1}^{r_{a}} 
\Big(1-q^{-d^{M'\cap M_{a,j}}(u_{j}')}\prod_{\alpha_{i}\in \Delta_{M_{a,j}}} t_{i}\Big),
\end{align*}
where $u'=(u'_{1},\cdots, u'_{r_{a}})$ is the coordinate with respect to the decomposition $\km_{(a)}=\km_{a, 1}\oplus\cdots\oplus \km_{a, r_{a}}$. The bound on the denominator of $g_{A}^{G}(u;\btt)$ then follows, and with (\ref{ga reduce 1}) and (\ref{ga reduce 2}) we can bound the denominator of $g_{A}^{M}(u;\btt)$.

\end{proof}

\subsection{Generating series of the affine Springer fibers}

We resume the discussions in \S \ref{counting fund domain}. By theorem \ref{f to w}, we have
\begin{equation}\label{f to w 2}
|F_{\gamma}(\fq)|=|D(\gamma)|^{-\frac{1}{2}} 
\sum_{\substack{M,\,L\in \cl(A)\\ M\subset L}}
(-1)^{\dim(\ka_{A}^{M})}J_{A}^{M}(\gamma, \mathbbm{1}_{\km\cap \kk})\,\rv_{M}^{L}(\ec_{A}(x_{0}))\, \re_{L},
\end{equation}
By proposition \ref{gkmreg}, $\ec_{A}(x_{0})$ is determined by the root valuation datum $\bnn$ of $\gamma$. To emphasize its dependence on the root valuation datum, we denote it by $\Pi_{\bnn}$.

\begin{thm}\label{fundamental domain}

Let $\bg=\gl_{d+1}$, let $\gamma$ be a regular diagonal matrix, then the number of points $|F_{\gamma}(\fq)|$ depends only on the root valuation datum $\bnn$ of $\gamma$, denoted  $F_{\bnn}(q)$. The generating series
$$
F(\mathbf{t}):=\sum_{n_{1}=0}^{+\infty}\cdots \sum_{n_{d}=0}^{+\infty}  F_{(n_{1},\cdots,n_{d})}(q) t_{1}^{n_{1}}\cdots t_{d}^{n_{d}}\in \bz [\![ t_{1},\cdots,t_{d}]\!]
$$
is a rational fraction, and its denominator divides $\prod_{i=1}^{d}(1-t_{i})$ times a certain power non exceeding $d$ of
$$
\prod_{\substack{L\in \cl(A)\\ L=L_{1}\times\cdots\times L_{r}}} \prod_{M\in \cl^{L}(A)}\prod_{\substack{u\in [\nn_{M}(F)] \\u=(u_{1},\cdots, u_{r})}}
\prod_{j=1}^{r} 
\Big(1-q^{d^{L_{j}}(0)-d^{M\cap L_{j}}(u_{j})}\prod_{\alpha_{i}\in \Delta_{L_{j}}} t_{i}\Big),
$$
where $L=L_{1}\times\cdots \times L_{r}$ is the factorization into linear groups, and $u=(u_{1},\cdots, u_{r})$ is the coordinate with respect to the corresponding decomposition of $\kl$. 

\end{thm}

\begin{proof}

It is clear that $|D(\gamma)|^{-\frac{1}{2}}$ depends only on the root valuation datum $\bnn$ of $\gamma$, we denote it by $D_{\bnn}$. By theorem \ref{main weighted}, the weighted orbital integrals $J_{A}^{M}(\gamma, \mathbbm{1}_{\km\cap \kk})$ depend only on $\bnn$ as well. By equation (\ref{f to w 2}), the counting point $|F_{\gamma}(\fq)|$ then depends only on $\bnn$. For the second assertion, it is enough to show that the generating series for $D_{\bnn}\cdot J_{A}^{M}(\bnn, \mathbbm{1}_{\km\cap \kk})\cdot\rv_{M}^{L}(\Pi_{\bnn})$ is rational. By equation (\ref{proj vol}), we have
$$
\rv_{M}^{L}(\Pi_{\bnn})=\sum_{M'\in \cl^{L}(A)}
\theta_{A}^{L}(M,M') \rv_{A}^{M'}(\Pi_{\bnn}),
$$
with $\theta_{A}^{L}(M,M')\neq 0$ if and only if the natural projection
$$
\ka_{A}^{M}\oplus \ka_{A}^{M'}\to \ka_{A}^{L}
$$
is an isomorphism. We call $M,M'$ \emph{transversal} if the condition holds for some $L\in \cl(A)$. The question is then reduced to showing that the generating series for 
$D_{\bnn}\cdot J_{A}^{M}(\bnn, \mathbbm{1}_{\km\cap \kk})\cdot\rv_{A}^{M'}(\Pi_{\bnn})$
is rational for all transversal  $M,M'$. With the Arthur-Shalika germ expansion, the question is further reduced to the rationality of the generating series of $D_{\bnn}\cdot  g_{A}^{M''}(\bnn, u)\cdot \rv_{A}^{M'}(\Pi_{\bnn})$ for $M''\in \cl^{M}(A)$ and $u\in [\nn_{M''}(F)]$. Notice that $M',M''$ are transversal again, and that if $\gamma$ is in minimal form, then it is also in minimal form for the groups $M',M''$ with respect to $\mathbf{B}_{0}\cap M'$ and $\mathbf{B}_{0}\cap M''$ respectively. 

To show the rationality, we look at the recursion of $D_{\bnn}\cdot  g_{A}^{M''}(\bnn, u)\cdot \rv_{A}^{M'}(\Pi_{\bnn})$ with respect to the transition $\bnn\mapsto \bnn+\mathbf{1}$.
For the factor $D_{\bnn}$, it is clear that
$$
D_{\bnn+\mathbf{1}}=q^{|\Phi^{+}(G,A)|}D_{\bnn}=q^{d_{G}(0)}D_{\bnn}.
$$
For the factor $\rv_{A}^{M'}(\Pi_{\bnn})$, notice that the $(G,A)$-family associated to $-\ec(x_{0})$ is equal to $\rr(\lambda;\,0,\gamma)$ up to multiplication by a global function on $(\ka_{A}^{G})^{*}$. Then similar to \cite{a1}, lemma 10.1, with the product formula (\ref{arthur prod simp}), we can show
\begin{equation*}
\rv_{A}^{M'}(\Pi_{\bnn+{\mathbf{1}}})=\sum_{J\in \cl^{M'}(A)}\rv_{A}^{J}(\Pi_{\bnn}) {\rc}_{J}^{M'}(0,\ep),
\end{equation*}
as $\Pi_{\bnn}$ satisfies the extra conditions on the facets leading to equation (\ref{arthur prod simp}). 
By theorem \ref{shalika homogeneity}, we get the recursion of $g_{A}^{M''}(\bnn, u)$ for $\bnn\mapsto \bnn+\mathbf{1}$. Combining them, we get
\begin{align}
&D_{\bnn+\mathbf{1}}\cdot \rv_{A}^{M'}(\Pi_{\bnn+\mathbf{1}})\cdot g_{A}^{M''}(\bnn+\mathbf{1},u) 
=q^{d_{G}(0)-d^{M''}(u)}\cdot \nonumber
\\
&
\sum_{J'\in \cl^{M'}(A)}\sum_{J''\in \cl^{M''}(A)}\sum_{v\in [\nn_{J''}(F)]}
D_{\bnn}\cdot \rv_{A}^{J'}(\Pi_{\bnn}) \cdot g_{A}^{J''}(\bnn,v) {\rc}_{J'}^{M'}(0,\ep) {\rc}_{J''}^{M''}(v,\ep)[v_{J''}^{M''}:u/\ep]. \label{recursion f}
\end{align}
Note that $J'$ and $J''$ are transversal, like $M'$ and $M''$. 
Following the same strategy of theorem \ref{shalika generating series}, we introduce 
$$
\bs^{M', M''}_{\Phi_{\bullet}}(u; \mathbf{t}):=\sum_{{\bnn\in \bz_{\ge 0}^{d},\, \Phi_{\bullet}(\bnn)=\Phi_{\bullet}}}  
D_{\bnn}\cdot \rv_{A}^{M'}(\Pi_{\bnn})\cdot g_{A}^{M''}(\bnn,u) t_{1}^{n_{1}}\cdots t_{d}^{n_{d}},
$$
and we break the summation in $\bs^{M',M''}_{\Phi_{\bullet}}(u;\mathbf{t})$ into two terms, $\bs_{1}^{M', M''}(u;\mathbf{t})$ and $\bs_{2}^{M', M''}(u;\mathbf{t})$, according to whether $m_{1}=0$ or not. 
Let $l$ be the length of the filtration $\Phi_{\bullet}$. For $i=1, \cdots, l$, let $M_{(i)}$ be the Levi subgroup with root system $\Phi_{i}$.

For the term $\bs_{1}^{M', M''}(u;\mathbf{t})$, 
suppose that $M_{(2)}$ has a factorization into linear groups $M_{(2)}=M_{2, 1}\times\cdots M_{2,r_{2}}$. Let $A_{i}=A\cap M_{2, i}$ and $\bnn_{i}$ be the restriction of the root valuation datum $\bnn$ to the simple roots in $\Delta_{M_{2,i}}$. As $m_{1}=0$, we have
\begin{equation}\label{factorize volume}
D_{\bnn}=\prod_{i=1}^{r_{2}}D^{M_{2,i}}_{\bnn_{i}}\quad \text{and}\quad\rv_{A}^{M'}(\Pi_{\bnn})=\prod_{i=1}^{r_{2}}\rv_{A_{i}}^{M'\cap M_{2, i}}(\Pi_{\bnn_{i}}).
\end{equation}
By lemma \ref{reduce gm levi 1}, there is a linear expression of $g_{A}^{M''}(\bnn, u)$ in terms of 
$$
g_{A}^{M'''}(\bnn', v'), \quad M'''\in \cl^{M''\cap M_{(2)}}(A), v'\in [\nn_{M'''}(F)],
$$ 
where $\bnn'$ is the restriction of the root valuation datum $\bnn$ to the simple roots in $\Delta_{M_{(2)}}$.
Let $M'''_{i}=M'''\cap M_{2, i}$, and let $v'=(v'_{1},\cdots, v'_{r_{2}})$ with respect to the decomposition $\km_{(2)}=\km_{2,1}\oplus\cdots\oplus \km_{2,r_{2}}$. With the factorization (\ref{factorize g}), we have 
\begin{equation}\label{factorize gm}
g_{A}^{M'''}(\bnn', v')=\prod_{i=1}^{r_{2}}g_{A_{i}}^{M'''\cap M_{2, i}}(\bnn_{i}, v'_{i}).
\end{equation}
With the equation (\ref{factorize volume}) and (\ref{factorize gm}), we obtain that

\begin{lem}\label{rationality S1}

The power series $\bs_{1}^{M',M''}(u;\btt)$ can be written as a linear combination of 
$$
\prod_{i=1}^{r_{2}} \bs_{\Phi_{\bullet}^{M_{2,i}}}^{M'\cap M_{2, i}, M'''\cap M_{2,i}}(v_{i}', \btt_{i}), \quad M'''\in \cl^{M''\cap M_{(2)}}(A), v'\in [\nn_{M'''}(F)],
$$
where $\Phi_{\bullet}^{M_{2,i}}$ is the reduced strict filtration on $\Phi(M_{2, i}, A_{i})$ induced by $\Phi_{\bullet}$, and $\btt_{i}$ is indeterminants indexed by $\Delta_{M_{2, i}}$.

\end{lem}

\noindent The rationality of $\bs_{1}^{M', M''}(u,\btt)$ then follows by induction on the rank of $\bfg$.

We proceed to the rationality of $\bs_{2}^{M', M''}(u, \btt)$. With the same lines of reasoning as in the proof of theorem \ref{shalika generating series} and with the recursion (\ref{recursion f}), we have

\begin{align*}
&\bs^{M', M''}_{\Phi_{\bullet}}(u;\mathbf{t})
\\
&=\bs^{M', M''}_{1}(u;\mathbf{t})+\sum_{{\bnn\in \bz_{\ge 0}^{d},\, \Phi_{\bullet}(\bnn)=\Phi_{\bullet}}}  
D_{\bnn+\mathbf{1}}\cdot \rv_{A}^{M'}(\Pi_{\bnn+\mathbf{1}})\cdot g_{A}^{M''}(\bnn+\mathbf{1},u)  t_{1}^{n_{1}+1}\cdots t_{d}^{n_{d}+1}
\\
&=\bs^{M', M''}_{1}(u;\mathbf{t})+\sum_{{\bnn\in \bz_{\ge 0}^{d},\, \Phi_{\bullet}(\bnn)=\Phi_{\bullet}}}  
q^{d^{G}(0)-d^{M''}(u)}\sum_{J'\in \cl^{M'}(A)}\sum_{J''\in \cl^{M''}(A)}\sum_{v\in [\nn_{J''}(F)]}
\\
&\hspace{3cm}D_{\bnn}\cdot \rv_{A}^{J'}(\Pi_{\bnn}) \cdot g_{A}^{J''}(\bnn,v) \nonumber
 \cdot {\rc}_{J'}^{M'}(0,\ep) {\rc}_{J''}^{M''}(v,\ep)[v_{J''}^{M''}:u/\ep]
t_{1}^{n_{1}+1}\cdots t_{d}^{n_{d}+1}
\\
&=\bs^{M', M''}_{1}(u;\mathbf{t})+  
q^{d^{G}(0)-d^{M''}(u)} \Bigg[\bs_{\Phi_{\bullet}}^{M', M''}(u/\ep;\mathbf{t})t_{1}\cdots t_{d}+ 
\sum_{\substack{J'\in \cl^{M'}(A), \,J''\in \cl^{M''}(A)\\ (J', \,J'')\neq (M',\, M'')}}\sum_{v\in [\nn_{J''}(F)]} \\
&\hspace{6cm}  
\bs_{\Phi_{\bullet}}^{J', J''}(v, \btt)
\cdot {\rc}_{J'}^{M'}(0,\ep) {\rc}_{J''}^{M''}(v,\ep)[v_{J''}^{M''}:u/\ep]
t_{1}\cdots t_{d}\Bigg].
\end{align*}
Since $u$ and $u/\ep$ belongs to the same nilpotent orbit in $\nn_{M''}(F)$, we have 
$$
\bs^{M', M''}_{\Phi_{\bullet}}(u;\mathbf{t})=\bs^{M', M''}_{\Phi_{\bullet}}(u/\ep;\mathbf{t}),
$$
and the above equation reduces to
\begin{align}
&\Big(1-q^{d^{G}(0)-d^{M''}(u)}t_{1}\cdots t_{d} \Big)\bs^{M', M''}_{\Phi_{\bullet}}(u;\mathbf{t}) \nonumber
\\
&=\bs^{M', M''}_{1}(u;\mathbf{t})+  
q^{d^{G}(0)-d^{M''}(u)}
\sum_{\substack{J'\in \cl^{M'}(A), \,J''\in \cl^{M''}(A)\\ (J', \,J'')\neq (M',\, M'')}}\sum_{v\in [\nn_{J''}(F)]} \nonumber
\\
&\hspace{5cm}  
\bs_{\Phi_{\bullet}}^{J', J''}(v, \btt)
\cdot {\rc}_{J'}^{M'}(0,\ep) {\rc}_{J''}^{M''}(v,\ep)[v_{J''}^{M''}:u/\ep]
t_{1}\cdots t_{d}. \label{recursion f sum}
\end{align}
The rationality of $\bs^{M', M''}_{\Phi_{\bullet}}(u;\mathbf{t})$ then follows by induction on the rank of both $M'$ and $M''$.  

To bound the denominator of $\bs^{M', M''}_{\Phi_{\bullet}}(u;\mathbf{t})$, with the recursion (\ref{recursion f sum}), we obtain that it divides $\big(1-q^{d^{G}(0)-d^{M''}(u)}t_{1}\cdots t_{d} \big)$ times those of $\bs_{1}^{M',M''}(u; \btt)$ and $\bs_{\Phi_{\bullet}}^{J', J''}(v, \btt)$. By lemma \ref{rationality S1}, the denominator of $\bs_{1}^{M',M''}(u; \btt)$ can be further reduced. At the end, we arrive at $\bs_{\Phi_{\bullet}}^{A, A}(0, \btt)$, which has denominator $\prod_{i=1}^{d}(1-t_{i})$. For $a=1,\cdots, l$, let $M_{(a)}=M_{a, 1}\times\cdots \times M_{a, r_{a}}$ be the factorization into linear groups, we conclude that the denominator of $\bs^{M', M''}_{\Phi_{\bullet}}(u;\mathbf{t})$ divides $\prod_{i=1}^{d}(1-t_{i})$ times a certain power non exceeding $d$ of
\begin{align*}
&\prod_{J''\in \cl^{M''}(A)}\prod_{v\in [\nn_{J''}(F)]} 
\Big(1-q^{d^{G}(0)-d^{J''}(v)}t_{1}\cdots t_{d} \Big)\cdot
\\
&\prod_{a=1}^{l} \prod_{M'''\in \cl^{M''\cap M_{(a)}}(A)}\prod_{\substack{u'\in [\nn_{M'''}(F)]\\ u'=(u'_{1},\cdots, u'_{r_{a}})}}
\prod_{j=1}^{r_{a}} 
\Big(1-q^{d^{M_{a,j}}(0)-d^{M'''\cap M_{a,j}}(u_{j}')}\prod_{\alpha_{i}\in \Delta_{M_{a,j}}} t_{i}\Big),\end{align*}
where $u'=(u'_{1},\cdots, u'_{r_{a}})$ is the coordinate with respect to $\km_{(a)}=\km_{a, 1}\oplus\cdots\oplus \km_{a, r_{a}}$, and the existence of a power non exceeding $d$ is due to the contributions of items indexed by $(J', M''), J'\neq M'$ at the right hand side of the recursion (\ref{recursion f sum}). The theorem then follows by summing $\bs^{M', M''}_{\Phi_{\bullet}}(u;\mathbf{t})$ together.

\end{proof}

\begin{rem}

For the goup $\bfg=\gl_{3}$ and the split elements $\gamma\in \ka$, we have constructed the affine pavings of the fundamental domains $F_{\gamma}$, and calculate their Poincar\'e polynomials in our work \cite{chen1}. As a corollary, we get 
\begin{align}
&\sum_{n_{1}=1}^{+\infty}\sum_{n_{2}=1}^{+\infty}  P_{(n_{1},n_{2})}(t)\, T_{1}^{n_{1}} T_{2}^{n_{2}}
\nonumber
\\
&=\lefteqn{2\Bigg\{\frac{(t^{2}+1)T_{1}T_{2}}{(1-T_{2})(1-T_{1}T_{2})(1-t^{4}T_{1}T_{2})^{2}}+\frac{t^{4}T_{1}T_{2}^{2}(3-t^{4}T_{1}T_{2})}{(1-T_{2})(1-t^{2}T_{2})(1-t^{4}T_{1}T_{2})^{2}}} \nonumber
\\
&+\frac{4t^{4}T_{1}T_{2}}{(1-t^{2}T_{2})(1-t^{4}T_{1}T_{2})^{2}(1-t^{6}T_{1}T_{2})}+\frac{t^{6}T_{1}T_{2}}{(1-t^{2}T_{2})(1-t^{6}T_{1}T_{2})}\Bigg\} \nonumber
\\
&-\left[\frac{(t^{2}+1)T_{1}T_{2}}{(1-T_{1}T_{2})(1-t^{4}T_{1}T_{2})^{2}}+\frac{4t^{4}T_{1}T_{2}}{(1-t^{4}T_{1}T_{2})^{2}(1-t^{6}T_{1}T_{2})}+\frac{t^{6}T_{1}T_{2}}{1-t^{6}T_{1}T_{2}}\right]. \label{P generating series}
\end{align}
As the fundamental domains admit affine pavings, their Poincar\'e polynomials match the counting point formula if we set $q=t^{2}$. Compare the denominator in the above formula with that of theorem \ref{fundamental domain}, we find that the former matches part of the items in the latter. Let $u_{0}$ and $u_{1}$ be the regular and subregular nilpotent matrix respectively, i.e.
$$
u_{0}=\begin{bmatrix}
0&1&\\
&0&1\\
&&0
\end{bmatrix},\quad
u_{1}=\begin{bmatrix}
0&&\\
&0&1\\
&&0
\end{bmatrix},
$$
then $\cn_{G}$ consists of $3$ orbits, that of $u_{0}, u_{1}$ and $0$. We calculate that $$d^{G}(u_{0})=0, \quad d^{G}(u_{1})=1, \quad d^{G}(0)=3.
$$ 
Let $M\in \cl(A)$ be the Levi subgroup with roots $\{\pm \alpha_{2}\}$, the nilpotent cone $\cn_{M}$ consists of two nilpotent orbits, that of $0$ and $u_{1}$. We calculate 
$$
d^{M}(0)=1, \quad d^{M}(u_{1})=0.
$$ 
In the formula of theorem \ref{fundamental domain}, we find denominators
\begin{align*}
&\big(1-q^{d^{G}(0)-d^{G}(u_{0})}t_{1}t_{2}\big)\big(1-q^{d^{G}(0)-d^{G}(u_{1})}t_{1}t_{2}\big)\big(1-q^{d^{G}(0)-d^{G}(0)}t_{1}t_{2}\big)
\\
&\cdot (1-q^{d^{M}(0)-d^{M}(u_{1})}t_{2})(1-q^{d^{M}(0)-d^{M}(u_{0})}t_{2})\\
&=(1-q^{3}t_{1}t_{2}) (1-q^{2}t_{1}t_{2})(1-t_{1}t_{2})\cdot (1-qt_{2})(1-t_{2}).
\end{align*}
They actually appear in the formula (\ref{P generating series}) if we set $q=t^{2}$.

\end{rem}

\begin{cor}\label{conj rational reduced}

The conjecture $\ref{conj rational}$ holds, if we assume the hypothesis that the affine Springer fibers for the split elements are cohomologically pure. 

\end{cor}

\begin{proof}

By theorem 1.2 of our work \cite{chen1}, if the affine Springer fibers are cohomologically pure,  then their fundamental domains will be cohomologically pure as well. For $\gamma\in \ka$, the maximal torus $A$ acts on $F_{\gamma}$ and the $A$-equivariant homology of $F_{\gamma}$ can be calculated with the work \cite{cl1} of Chaudouard-Laumon. Let 
\begin{align*}
\cd&=H^{\bullet}_{A}(\mathrm{pt}, \mathbf{Q}_{\ell})=\bigoplus_{n\in \bn_{0}}\mathrm{Sym}^{n}(X^{*}(A), \mathbf{Q}_{\ell})(-n),
\\
\cs&=H^{A}_{\bullet}(\mathrm{pt}, \mathbf{Q}_{\ell})=\bigoplus_{n\in \bn_{0}}\mathrm{Sym}^{n}(X_{*}(A), \mathbf{Q}_{\ell})(n).
\end{align*}
With the cap product, the homology groups $\cs$ is naturally a graded $\cd$-module.
Forget the Tate twist, we can identify $\cs$ as the ring of polynomial functions on $\ka^{*}$ and $\cd$ as the ring of linear differential operators on $\ka^{*}$. For $\alpha\in \Phi(G,A)$, let $\partial_{\alpha}$ be the corresponding differential operator on $\ka^{*}$. 
The set of fixed points $\xx_{\gamma}^{A}$ can be naturally identified with $X_{*}(A)$, hence $F_{\gamma}^{A}$ can be identified as a subset of $X_{*}(A)$. For $\alpha\in \Phi^{+}(G,A)$, let 
$$
R_{\alpha, i}=\sum_{\lambda \text{ satisfying } (*)} (1-\alpha^{\vee})^{i}\lambda\otimes \cs\{\partial_{\alpha}^{i}\}\subset F_{\gamma}^{A}\otimes \cs,
$$
where $(*)$ refers to the condition that $\lambda, \alpha^{\vee}\lambda, \cdots, (\alpha^{\vee})^{i}\lambda\in F_{\gamma}^{A}$. 
By proposition 10.3 of \cite{cl1}, the $A$-equivariant homology of $F_{\gamma}$ fits into the exact sequence
$$
0\to \sum_{\alpha\in \Phi^{+}(G, A)}\sum_{i=1}^{\val(\alpha(\gamma))}R_{\alpha, i}\to F_{\gamma}^{A}\otimes \cs_{\bullet}\to H^{A}_{\bullet}(F_{\gamma}, \mathbf{Q}_{\ell})\to 0.
$$
Moreover, the cohomological purity of $F_{\gamma}$ implies that
$$
H_{\bullet}(F_{\gamma}, \mathbf{Q}_{\ell})=H^{A}_{\bullet}(F_{\gamma}, \mathbf{Q}_{\ell})\{I\},
$$
where $I\subset \cd$ is the augmentation ideal. We conclude that $H^{*}(F_{\gamma}, \mathbf{Q}_{\ell})$ is non-trivial only in the even degree and the Frobenius action on $H^{2i}(F_{\gamma}, \mathbf{Q}_{\ell})$ is scalar product by $q^{i}$. The Poincar\'e polynomial $P(F_{\gamma},t)$ of $F_{\gamma}$ is then a polynomial in $t^{2}$, and with the fixed point formula of Grothendieck-Lefschetz we deduce that 
$$
|F_{\gamma}(\fq)|=P(F_{\gamma}, q^{1/2}).
$$
The conjecture \ref{conj rational} then follows from theorem \ref{fundamental domain}.

\end{proof}

\begin{rem}
It is questionable whether the methods of this section work for general semisimple elements. 
Let $\gamma=\diag(\gamma_{1},\cdots, \gamma_{r})$ be a regular semisimple element such that $\gamma_{i}$ is anisotropic and non-split for all $i$. 
Let $M_{0}$ be the minimal Levi subgroup containing $\gamma$. 
For $u\in [\cn_{G}(F)]$, with the homogeneity of the Arthur-Shalika germ, we can reduce $g_{M_{0}}^{G}(\gamma, u)$ to $g_{M_{0}}^{M}(\ep^{-n}\gamma, v), M\in \cl(M_{0}), v\in [\cn_{M}(F)]$ for a certain $n\in \bn_{0}$ such that the smallest root valuation of $\ep^{-n}\gamma$ lies in the interval $(0,1)$.
The problem is that the parabolic descent for the Arthur-Shalika germ, or more precisely the theorem \ref{germ descent}, doesn't work for such elements. Indeed, $\ep^{-n}\gamma$ doesn't lie in any invariant open neighbourhood $U_{s}^{+}, s\in \ka$.
For the generalization, we will need a germ expansion around any semisimple element, but this is much harder.

\end{rem}

\bigskip
{\small
\noindent\textbf{Acknowledgement:}
We want to thank Prof. G\'erard Laumon and Prof. Jean-Loup Waldspurger for their interest in the work. We are very grateful to Prof. Waldspurger for pointing out the errors in the preliminary versions of this paper and  for his very kind suggestion for the proof of lemma \ref{key lemma}. 
We also want to thank an anonymous referee for the suggestion that it should be possible to bound the denominators of the generating series with the recursion relations. 
This work is partially supported by Tsinghua initiative research funding.

}

\bigskip
\small
\noindent
\begin{tabular}{lll}
&Zongbin {\sc Chen} &\\ 
\\
&Current address: &Previous address:\\
&&\\
&School of mathematics, &Yau Mathematical Science Center,\\
&Shandong University,&Tsinghua University,\\
&Shanda south road 27, Licheng district, JiNan, 250100, &Haidian, 100084,\\
&Shandong, China &Beijing, China \\
&email: {\tt zongbin.chen@email.sdu.edu.cn} &

\end{tabular}


\begin{thebibliography}{100}
\labelwidth=4em
\addtolength\leftskip{20pt}
\setlength\labelsep{0pt}



\bibitem[A1]{a}{J. Arthur, \textit{The characters of discrete series as orbital integrals}, Invent. math. \textbf{32}(1976), 205-261.}  

\bibitem[A2]{a4}{J. Arthur, \textit{The trace formula in invariant form}, Ann. of Math., 114 (1981), 1-74.}






\bibitem[A3]{a1}{J. Arthur, \textit{The local behavior of weighted orbital integrals}, Duke. Math. J., Vol 56, No. 2, 1988.}

\bibitem[A4]{a2}{J. Arthur, \textit{An introduction to the trace formula}. Harmonic analysis, the trace formula, and Shimura varieties, 1–263, Clay Math. Proc., 4, Amer. Math. Soc., Providence, RI, 2005.} 

\bibitem[B]{bourbaki lie}{N. Bourbaki, \textit{Groupes et alg\`ebres de Lie}. Chapitre 4-6.  Masson, Paris, 1982.} 



\bibitem[Bez]{be}{R. Bezrukavnikov, \textit{The dimension of the fixed point set on affine flag manifolds}, Math. Res. Lett. 3 (1996), 185–189.}
















\bibitem[C1]{chen1}{Z. Chen, \textit{On the fundamental domain of affine Springer fibers}. Math. Z. 286 (2017), no. 3-4, 1323-1356.}


\bibitem[C2]{chen2}{Z. Chen, \textit{Truncated affine Springer fibers and Arthur's weighted orbital integrals}, J. Inst. Math. Jussieu 22 (2023), no. 4, 1757-1818.}

\bibitem[CL1]{cl1}{P-H. Chaudouard, G. Laumon, \textit{Sur l'homologie des fibres de Springer affines tronqu\'ees}, Duke Math. J. 145 (2008), no. 3, 443–535.}


\bibitem[CL2]{cl2}{P-H. Chaudouard, G. Laumon, \textit{Le lemme fondamental pond\'er\'e. I. Constructions g\'eom\'etriques.} Compos. Math. 146 (2010), no. 6, 1416–1506.} 

\bibitem[CM]{nilpotent ref}{D. Collingwood, W. McGovern, 
Nilpotent orbits in semisimple Lie algebras.
Van Nostrand Reinhold Math. Ser.
Van Nostrand Reinhold Co., New York, 1993.}


\bibitem[D1]{debacker gl}{S. DeBacker, \textit{Homogeneity of certain invariant distributions on the Lie algebra of $p$-adic $\gl_{n}$}, Compos. Math. 124 (1) (2000) 11-16.}

\bibitem[D2]{debacker}{S. DeBacker, \textit{Homogeneity results for invariant distributions of a reductive $p$-adic group.} Ann. Sci. \'Ecole Norm. Sup. (4) 35 (2002), no. 3, 391-422.}




\bibitem[GKM1]{gkm1}{M. Goresky, R. Kottwitz, R. MacPherson, \textit{Homology of affine Springer fibers in the unramified case}, Duke Math. J. 121 (2004), no. 3, 509-561.}


\bibitem[GKM2]{gkm3}{M. Goresky, R. Kottwitz, R. MacPherson, \textit{Regular points in affine Springer fibers}, Michigan Math. J. 53 (2005), no. 1, 97-107.}

\bibitem[GKM3]{gkm4}{M. Goresky, R. Kottwitz, R. MacPherson, \textit{Codimensions of root valuation strata}, Pure. Appl. Math. Q. 5 (2009), no. 4, 1253-1310.}








\bibitem[K1]{k1}{R. E. Kottwitz, Isocrystals with additional structure. Compositio Math. 56 (1985), no. 2, 201–220.}




\bibitem[K2]{kottwitz course}{R. E. Kottwitz, \textit{Harmonic analysis on reductive $p$-adic groups and Lie algebras}. Harmonic analysis, the trace formula, and Shimura varieties, 393–522, Clay Math. Proc., 4, Amer. Math. Soc., Providence, RI, 2005}


\bibitem[KL]{kl}{D. Kazhdan, G. Lusztig, \textit{Fixed point varieties on affine flag manifolds}, Israel. J. Math. \textbf{62}(1988), 129-168.}



\bibitem[KV]{kv}{R. Kottwitz, E. Viehman, \textit{Generalized affine Springer fibers}, J. Inst. Math. Jussieu 11 (2012), no. 3, 569–609.}






\bibitem[MP]{mp}{A. Moy, G. Prasad, \textit{Unrefined minimal K-types for p-adic groups}, Invent. Math. 116 (1994), 393–408.}











\bibitem[W]{walds local trace}{J.-L. Waldspurger, \textit{Une formule des traces locale pour les alg\`ebres de Lie p-adiques}. J. Reine Angew. Math. 465 (1995), 41–99.}

\bibitem[We]{weil integration}{A. Weil, Adeles and algebraic groups, with appendices by M. Demazure and Takashi Ono. Progress in Mathematics, 23. Birkh\"auser, Boston, Mass., 1982.}


\end{thebibliography}
\end{document}